\documentclass{article}
\usepackage{color}
\usepackage{amssymb}
\usepackage{amstext}

\newtheorem{Corollary}{Corollary}
\newtheorem{theorem}{Theorem}
\newtheorem{definition}{Definition}
\newtheorem{proposition}{Proposition}
\newtheorem{lemma}{Lemma}

\usepackage[dvips]{graphicx}
\newcommand{\K}{$\Bbb K$}
\newcommand{\C}{\mathbb C}
\newcommand{\g}{\frak{g}}
\newcommand{\w}{\frak{r}}
\newcommand{\R}{\mathbb R}
\newcommand{\Z}{\mathbb Z}

\newcommand{\n}{\frak{n}}
\newcommand{\m}{\frak{t}}
\newcommand{\s}{\frak{s}}

\begin{document}

{\noindent \Huge Lie algebras : \\
\medskip \\
Classification, Deformations and Rigidity}

\medskip

\noindent \rule[1cm]{14cm}{0.1cm}
\medskip

\noindent Michel GOZE, Universit\'e de Haute Alsace, MULHOUSE (France). M.Goze@uha.fr
\bigskip

\noindent Lessons given during the {\it Cinqui\`eme Ecole de G\'eom\'etrie Diff\'erentielle et
Syst\`emes Dynamiques} , ENSET ORAN (Algeria), november 4-11, 2006.

\medskip

\noindent \rule[1cm]{14cm}{0.01cm}

\bigskip

{\small
In the first section we recall some basic notions on  Lie algebras. In  a second time w
e study the algebraic variety of complex $n$-dimensional Lie algebras. We present different notions of deformations :
Gerstenhaber deformations, pertubations, valued deformations and we  use these tools to study 
some properties of this variety. Finaly we introduce the concept of rigidity and we present some results on the class
of rigid Lie algebras.
}

\medskip

\noindent{\bf Table of contents}

\noindent page 1 . Section 1 : Lie algebras.

page 1 : 1.1 Definitions and Examples.

page 2 : 1.2 The Lie algebra of a Lie group.

page 2 : 1.3 Lie admissible algebras.

\noindent page 3 . Section 2 : Classification of   Lie algebras.

page 4 : 2.1 Simple Lie algebras

page 5 : 2.2 Nilpotent Lie algebras

page 6 : Solvable Lie algebras

\noindent page 9 . Section 3 : The algebraic variety of complex Lie algebras

page 9 : 3.1 The algebraic variety $L_n$.

page 10 : 3.2 The tangent space to the orbit.

page 10 : 3.3 The tangent space to $L_n$.

\noindent page 11 . Section 4 : The scheme $\mathcal{L}_n$.

page 11 : 4.1  Definition

page 13 : 4.2  The tangent space to the scheme.

\noindent page 13 . Section 5 : Contractions of Lie algebras

page 13 : 5.1  Definition.

page 14 : 5.2 Examples.

page 16 : 5.3 In\"on\"u-Wigner contractions.

page 19 : 5.4 In\"on\"u-Wigner contractions of Lie groups

page 19 : 5.5 Wiemar-Woods contractions

page 19 : 5.6 The diagram of contractions

\bigskip

\begin{center}
\section{Lie algebras} 
\end{center}

\medskip

\subsection{ Definition and examples}

\medskip

In this lecture, the considered Lie algebras are complex (sometimes real but this case  will be precised).

\medskip

\framebox{\parbox{6.25in}{
\begin{definition}
A Lie algebra $\g$ is a pair $(V,\mu)$ where $V$ is a complex vector space and $\mu$ a bilinear map
$$
\mu :V \times V \rightarrow V
$$
satisfying :
$$
\mu (X,Y)=-\mu (Y,X),\quad \forall X,Y \in V,$$
$$
 \mu (X,\mu (Y,Z))+\mu (Y,\mu (Z,X))+\mu (Z,\mu (X,Y))=0,\quad \forall
X,Y,Z\in V.$$
\end{definition}
}}

\bigskip

\noindent This last identity is called the Jacobi relation.

\medskip 

\subsection{Examples}

\medskip

1. The simplest case is given taking $\mu=0$. Such a Lie algebra is called abelian. 

\medskip

\noindent 2. Suppose that $V$ is $2$-dimensional. Let $\{e_1,e_2\}$ be a basis of $V$. If we put
$$\mu(e_1,e_2)=ae_1+be_2$$
with the condition of skew-symmetry , this map satisfies the Jacobi condition and it is a multiplication of
Lie algebra.

\medskip

\noindent 3. Let $sl(2,\C)$ be the vector sapce of matrices $A$ of order $2$ such that $tr(A)=0$ where $tr$ indicates the
trace of the matrix $A$. The product
$$\mu(A,B)=AB-BA$$
is well defined on $sl(2,\C)$ because $tr(AB-BA)=0$ as soon as $A,B \in sl(2,\C)$. We can see that this product
satisfies the Jacobi condition and we have a Lie algebra of dimension $3$.

\subsection{The Lie algebra of a Lie group}

Let $G$ be a (complex) Lie group of dimension $n$. For every $g \in G$, we denote by $L_g$ the automorphism of $G$ given by
$$L_g(x)=gx.$$
It is called the left translation by $g$. Its differential map $(L_g)^*_x$ is an isomorphism
$$(L_g)^*_x: T_x(G)\longrightarrow T_{gx}(G)$$
where $T_x(G)$ designates the tangent space at $x$ to $G$. A vector field $X$ on $G$ is called left invariant if it satisfies
$$(L_g)^*_x(X(x))=X(gx)$$
for every $x$ and $g$ in $G$. We can prove that, if $[X,Y]$ denotes the classical braket of vector fields and if
$X$ and $Y$ are two left invariant vector fields of $G$ then $[X,Y]$ is also a left invariant vector field. This shows
that the vector space of left invariant vector fields of $G$ is provided with a Lie algebra structure denoted $L(G)$. 
As a left invariant vector field $X$ is defined as soon as we know its value $X(e)$ to the unity element $e$ of $G$ 
the vector space
$L(G)$ can be identified to $T_e(G)$. For example the Lie algebra of the algebraic group $SL(2,\C)$ is isomorphic to
$sl(2,\C)$. (We shall define the notion of isomorphism later).

\medskip

From this construction, to each Lie group we have only one associated Lie algebra. But the converse is not true. In fact
two Lie groups which are locally isomorphic have the same associated Lie algebra. On a other hand, we can construct from a 
finite dimensional complex Lie algebra a Lie group. This is a little bit complicated. But from the bracket of 
a finite dimensional Lie algebra $\g$, we can define a local group whose product is given by the Campbell-Hausdorff
formula:
$$X.Y= X+Y+1/2[X,Y]+1/12[[X,Y],Y]-1/12[[X,Y],X]+....$$
which is a infinite sequence of terms given from the bracket. This local structure can be extended to a global structure
of Lie group. Then we have a one-one correspondance between finite dimensional Lie algebra (on $\C$ for example or on $\R$)
and the simply connected, connected Lie groups. We can note that the dimension of the Lie group as differential manifold is
equal to the dimension of its Lie algebra as vector space.

\subsection{Relation between Lie group and its Lie agebra}

In other talks of this school, the notion of linear algebraic groups is often used. If the Lie group is a linear group, that
is a Lie subgroup of the Lie group $Gl(n,\C)$, then its Lie algebra is a Lie subalgebra of the Lie algebra
$gl(n,\C)$ whose elements are the complex matrices of order $n$.  In this case we can easily construct a map between the Lie
algebra and a corresponding Lie group. It is based on the classical exponential map.

The exponential map from the Lie algebra $gl(n,\C)$ of the Lie group $GL(n,\C)$  to $GL(n,\C)$  is defined 
by the usual power series:
$$Exp(A) = Id + A + \frac{A^2}{2!}+\frac{A^3}{3!}+  ...$$
for matrices $A$. If $G$ is any Lie subgroup of $GL(n,\C)$ then the exponential map takes the Lie algebra 
$\g$ of $G$ into $G$, so we have an exponential map for all matrix groups. 

The definition above is easy to use, but it is not defined for Lie groups that are not matrix groups. The Ado theorem
precises that every finite dimensional Lie algebra can be represented as a linear algebra. This means that , for a 
gievn $n$ dimensional Lie algeba $\g$, there exists
an integer $N$ such that the elements of $\g$ are written as matrices of order $N$. The problem comes because we
do not know $N$. Moreover , if we write $\g$ as a subalgebra of $gl(N,\C)$, we can exprim the exponential map, but this map
depend of the choosen representation of $\g$ as linear algebra. 
We can solve both problems using a more abstract definition of the exponential map that works for all Lie groups, 
as follows.  
Every vector $X$ in $\g$ determines a linear map from $\R$ to $\g$ taking $1$ to $X$, which can be 
thought of as a Lie algebra homomorphism. Since $\R$ is the Lie algebra of the simply connected Lie group $\R$,
 this induces a Lie group homomorphism 
$$c : \R \longrightarrow G$$
 such that 
$$c(s + t) = c(s)+c(t)$$
for all $s$ and $t$. The operation on the right hand side is the group multiplication in $G$.
 The formal similarity of this formula with the one valid for the exponential function justifies the definition
$$Exp(X) = c(1).$$
This is called the exponential map, and it maps the Lie algebra $\g$ into the Lie group $G$. 
It provides a diffeomorphism between a neighborhood  
of $0$ in $\g$ and a neighborhood of $Id$ in $G$. 
The exponential map from the Lie algebra to the Lie group is not always onto, even if the group is connected
(For example, the exponential map of $sl(2,\C)$ is not surjective.) But if $\g$ is nilpotent $Exp$ is bijective.

\medskip

\subsection{Lie admissible algebras}

Let $ \mathcal{A} $ be a complex associative algebra whose multiplication is denoted by  $A.B$ with $A, B \in 
\mathcal{A}$. It is easy to see that $$[A,B]=AB-BA$$
is a Lie bracket. Let $\mathcal{A}^L$ the corresponding Lie algebra. For example, if $M(n,\C)$ is the vector space
of complex matrices of order $n$, it is an associative algebra for the usual product of matrices then 
$[A,B]=AB-BA$ is a Lie product on $M(n,\C)$. This Lie algebra is usually denoted $gl(n,\C)$. We can note that there exists
Lie algebras which are not given by an associative algebra.

\medskip
\noindent
\framebox{\parbox{6.25in}{
\begin{definition}
A Lie-admissible algebra is a (nonassociative) algebra $\mathcal{A}$ whose product $A.B$ is such that
$$[A,B]=AB-BA$$
is a product of Lie algebra.  
\end{definition}
}}
 
\bigskip

It is equivalent to say that the product $A.B$ satisfies
$$
\begin{array}{l}
(A.B).C-A.(B.C)-(B.A).C+B.(A.C)-(A.C).B+A.(C.B)-(C.B).A+C.(B.A)\\
+(B.C).A-B.(C.A)+(C.A).B-C.(A.C)=0
\end{array}
$$
for all $A,B,C \in \mathcal{A}.$ An interesting example concerns the Pre-Lie algebras, that is (nonassociative) algebras
whose product satisfies
$$(A.B).C-A.(B.C)-(A.C).B+A.(C.B)=0.$$
This algebras, also called left-symmetric algebras, are studied to determinate on a Lie group the  flat left 
invariant affine connections with torsionfree. For example, abelian Pre-Lie algebras determinates the affine structures
(i.e. flat torsionfree left invariant affine connections on a Lie group) on the corresponding Lie algebra, that is the Lie
algebra given by the Pre-Lie algebra. But an abelian Pre-Lie algebra is associative commutative.  Let us give, for example
the commutative associative algebra of dimension $2$ given by
$$
e_1.e_1=e_1, \ e_1.e_2=e_2.e_1=e_2 , \ e_2.e_2=e_2.$$
The left translation $l_X$ of this multiplication is given for $X=ae_1+be_2$ by 
$$l_X(e_1)=ae_1+be_2, \ l_X(e_2)=(a+b)e_2.$$
The associated Lie algebra is abelian and admits a representation on the Lie algebra $Aff(\R^2)$ (here we are interested by the
real case) given by
$$
\left(
\begin{array}{ll}
l_X & X \\
0 & 0
\end{array}
\right)
\ =
\left(
\begin{array}{lll}
a & 0 & a \\
b & a+b & b \\
0 & 0 & 0
\end{array}
\right).
$$
The Lie group associated to this affine Lie algebra is found taking the exponential of this matrix. We obtain
$$
G=\{
\left(
\begin{array}{lll}
e^a & 0 & e^a-1 \\
e^a(e^b-1) & e^ae^b & e^a(e^b-1)\\
0 & 0 & 1
\end{array}
\right)
 \  a, b \in \R .\}
$$
This gives the following affine transformations
$$
\left\{
\begin{array}{l}
x \longrightarrow e^ax + e^a-1 \\
y \longrightarrow e^a(e^b-1)x+ e^ae^by + e^a(e^b-1)
\end{array}
\right.
$$

\medskip

\subsection{Infinite dimensional Lie algebras}

The theory of infinite dimensional Lie algebras, that is the underlying vector space is of infinite dimension, 
is a little bit different. For example there is stil no correspondance between Lie algebras and Lie groups. W.T Van Est
showed example of Banach Lie algebras which are not associated to infinite Lie groups. But some infinite Lie algebras 
play fundamental role. The Kac-Moody algebras are graded infinite Lie algebras defined by generators and relations and
which are constructed as finite dimensional simple Lie algebras. Another example is given by the Lie algebra of the vector
fields of a differentiable manifold $M$. The bracket is given by the Lie derivative. This Lie algebra is very complicated.
It is well studied in cases of $M=\R$ or $M=S^1$. In this case the Lie algebra is associated to the Lie group of diffeomorphisms
of $\R$ or $S^1$.
Another family of infinite Lie algebras are the Cartan Lie algebras. They are defined as Lie algebras of infinitesimal
diffeomorphisms (vector fileds) which leave invariant a structure as symplectic structure or contact structure. For example
let $(M,\Omega )$ a symplectic variety, that is $\Omega $ is a symplectic form on the differential manifold $M$. We consider
$$L(M,\Omega )=\{X \ {\mbox{\rm vector \ fields \ on}} \ M, L_X\Omega =0\}$$
where $$L_X\Omega =i(X)d\Omega +d(i(X)\Omega )=d(i(X)\Omega )=0$$
is the Lie derivate. Then  $L(M,\Omega)$ is a real infinite Lie algebra. It admits a subalgebra $L_0$
constituted of vector fields of $L(M,\Omega )$ with compact support. Then Lichnerowicz proved that every
finite dimensional Lie algebra of $L_0$ is reductive (that is a direct product of a semi simple sub algebra by an abelian center)
and every nonzero ideal is of infinite dimensional.
\medskip
\begin{center}
\section{Classifications of Lie algebras}
\end{center}

\medskip

\noindent
\framebox{\parbox{6.25in}{
\begin{definition}
Two Lie algebras $\g$ and $\g'$ of multiplication  $\mu $ and $\mu ^{\prime }$ $\in L^{n}$ are said isomorphic, if
there is $f\in Gl(n,\mathbb{C)}$ such that
$$
\mu ^{\prime }(X,Y)=f*\mu (X,Y)=f^{-1}(\mu (f(X),f(Y)))
$$
for all $X,Y\in \g.$ \end{definition}
}}

\bigskip

For example any $2$-dimensional Lie algebra is gven by
$$\mu(e_1,e_2)=ae_1+be_2.$$
Suppose $a$ or $b$ not zero, for example $b$, the change of basis
$$
\left\{
\begin{array}{l}
f(e_1)=1/b \ e_1 \\
f(e_2)= a/b e_1+ e_2
\end{array}
\right.
$$
defines the isomorphic law
$$\mu (e_1,e_2)=e_2.$$
We deduce that every two dimensional Lie algebra is abelian or isomorphic to the Lie algebra whose product is 
$$\mu (e_1,e_2)=e_2.$$

The classification of $n$ dimensional Lie algebras consits to describe a representative Lie algebra 
of each class of isomorphism.
Thus the classification of complex (or real) $2$-dimensional Lie algebras is given by the following result :

\medskip

\noindent
\framebox{\parbox{6.25in}{
\begin{proposition}
Every $2$-dimensional complex (or real)Lie algebra is isomorphic to one of the following:

- The abelian $2$-dimensional Lie algebra.

- The Lie algebra defined from $\mu (e_1,e_2)=e_2$.
\end{proposition}
}}

\bigskip

The general classification, that is the classification of Lie algebra of arbitrary dimension is a very complicated problem 
and it is today unsolved. We known this general classification only up the dimension $5$. Beyond this dimension 
we know only partial classifications. We are thus led to define particular classes of Lie algebras.

\medskip

\subsection{Simple Lie algebras}

\medskip

\noindent
\framebox{\parbox{6.25in}{
\begin{definition}
A Lie algebra $\g$ is called simple if it is of dimension greater or equal to $2$ and do not contain proper ideal, that
is ideal not trivial and not equal to $\g$.
\end{definition}
}}

\bigskip

\noindent The classification of simple complex Lie algebras is wellknown. It is rather old and it is due primarily to works of Elie Cartan, Dynkin or Killing.
To summarize this work, let us quote only the final result.

\noindent Every simple complex Lie algebra is 

 - either of classical type that is isomorphic to one of the following: $su(n,\C)$ (type $A_n$), $so(2n+1,\C)$ (type $B_n$),
$sp(n,\C)$ (type $C_n$), $so(2n,\C)$ (type $D_n$)

- or exceptional that is of type $E_6, \ E_7, \ E_8, \ F_4, \ G_2$.

\noindent One can read the definition  of these algebras for example in the book of J.P. Serre entitled SemiSimple Lie algebras.

\medskip

Another class of Lie algebras concerns the semi-simple Lie algebra. By definition, a semi-simple Lie algebra is a 
direct product of
simple Lie algebras. Briefly let us point out the notion of direct product. Let $\g_1$ and $\g_2$ two complex Lie algebras whose multiplications
are denoted by $\mu_1$ and $\mu_2$. The direct product $\g_1 \otimes \g_2$ is the Lie algebra whose underlying vectorial space is
$\g_1 \oplus \g_2$ and whose multiplication is 
$$\mu(X_1+X_2,Y_1+Y_2)=\mu_1(X_1,Y_1)+\mu_2(X_2,Y_2)$$
for any $X_1,Y_1 \in \g_1$ and $X_2,Y_2 \in \g_2$.
We deduce easely, from the classification of simple Lie algebras that from semi simple complex Lie algebras.

\subsection{Nilpotent  Lie algebras}

Let $\g$ be a Lie algebra. Let us consider the following sequence of ideals
$$
\left\{
\begin{array}{l}
{\mathcal{C}}^1(\g)=\mu(\g,\g) \\
{\mathcal{C}}^p(\g)=\mu({\mathcal{C}}^{p-1}(\g),\g) \ for \  p>2.
\end{array}
\right.
$$
where $\mu({\mathcal{C}}^{p-1}(\g),\g)$ is the subalgebra generated by the product $\mu(X,Y)$ with $X \in
{\mathcal{C}}^{p-1}(\g)$ and $Y \in \g$. This sequence satisfies
$${\mathcal{C}}^p(\g) \subset {\mathcal{C}}^{p-1}(\g), \ p >0$$
where ${\mathcal{C}}^0(\g)=\g.$

\medskip

\noindent
\framebox{\parbox{6.25in}{
\begin{definition}
A Lie algebra $\g$ is called nilpotent if there is $k$ such that
$${\mathcal{C}}^k(\g)=\{0\}.$$
If a such integer exists, he smallest  $k$ is called the index of nilpotency or nilindex of $\g$.
\end{definition}
}}

\medskip

\noindent {\bf Examples.}

1. The abelian algebra satisfies ${\mathcal{C}}^1(\g)=\{0\}$. It is nilpotent of nilindex $1$.

2. The $(2p+1)$-dimensional Heisenberg algebra is given in a basis $\{e_1,...,e_{2p+1}\}$ from
$$\mu(e_{2i+1},e_{2i+2})=e_{2p+1}$$
for $i=0,...,p-1$. This Lie algebra is nilpotent of nilindex $2$. In this case ${\mathcal{C}}^1(\g)$ is generated
by the center $\{e_{2p+1}\}$. Generally speaking, we call $2$-step nilpotent Lie algebra a nilpotent Lie algebra whose nilindex is equal to $2$.

3. For any nilpotent $n$-dimensional Lie algebra, the nilindex is bounded by $n-1$.

\medskip

\noindent
\framebox{\parbox{6.25in}{
\begin{definition}
A nilpotent Lie algebra is called filiform if its nilindex is equal to $n-1$.
\end{definition}
}}

\medskip

\noindent For example, the following $4$-dimensional Lie algebra given by
$$
\left\{
\begin{array}{l}
\mu(e_1,e_2)=e_3 \\
\mu(e_1,e_3)=e_4
\end{array}
\right.
$$
is filiform. 

\medskip

The classification of complex (and real) nilpotent Lie algebras is known up the dimension $7$. For the 
dimensions $3$,$4$,$5$ and $6$, there exists 
only a finite number of classes of isomorphism. For the dimensions $7$ and beyond, there exist an infinity of isomorphism classes. For example , for dimension $7$, we have $6$ families of one parameter of 
non isomorphic Lie algebras.

The classification of nilpotent Lie algebras of dimension less or equal to $7$ is given in the web site

\medskip

\noindent
\framebox{\parbox{6.25in}{
\begin{center}
http:/ /  www.math.uha.fr / \~ \  algebre
\end{center}
}}

\medskip

\subsection{Solvable  Lie algebras}

Let $\g$ be a Lie algebra. Let us consider the following sequence of ideals
$$
\left\{
\begin{array}{l}
{\mathcal{D}}^1(\g)=\mu(\g,\g) \\
{\mathcal{D}}^p(\g)=\mu({\mathcal{D}}^{p-1}(\g),{\mathcal{D}}^{p-1}(\g)) \ for \  p>2.
\end{array}
\right.
$$
 This sequence satisfies
$${\mathcal{D}}^p(\g) \subset {\mathcal{D}}^{p-1}(\g), \ p >0$$
where ${\mathcal{D}}^0(\g)=\g.$

\medskip

\noindent
\framebox{\parbox{6.25in}{
\begin{definition}
A Lie algebra $\g$ is called solvable if there exists $k$ such that
$${\mathcal{D}}^k(\g)=\{0\}.$$
If a such integer exists, the smallest  $k$ is called the index of solvability or solvindex of $\g$.
\end{definition}
}}

\medskip

\noindent {\bf Examples.}

1. Any nilpotent Lie algebra is solvable because
$${\mathcal{D}}^{p}(\g) \subset {\mathcal{C}}^{p}(\g).$$

2. The non abelian $2$-dimensional Lie algebra is solvable. 

3. We can construct solvable Lie algebra starting of a nilpotent Lie algebra. First we need to some definitions.

\medskip

\noindent
\framebox{\parbox{6.25in}{
\begin{definition}
A derivation of a Lie algebra $\g$ is a linear endomorphism $f$ satisfying
$$\mu(f(X),Y)+\mu(X,f(Y))=f(\mu(X,Y))$$
for every $X,Y \in \g$.
\end{definition}
}}

\medskip

\noindent For example the endomorphisms $adX$ given by
$$adX(Y)=\mu(X,Y)$$
are derivations (called inner derivations). Let us note that any inner derivation is singular because $X \in Ker(adX)$. 
A Lie algebra provided with a regular derivation is nilpotent. Let us note too that if $\g$ is a simple or 
semi simple Lie algebra, any derivation is inner.

Let $\g$ be a $n$-dimensional Lie algebra and $f$ a derivation not inner. Consider the vector space $\g'=\g \oplus \C$
of dimension $n+1$ and we denote by $e_{n+1}$ a basis of the complementary space $\C$. We define a multiplication on $\g'$ by
$$
\left\{
\begin{array}{l}
\mu'(X,Y)=\mu(X,Y), \ X,Y \in \g\\
\mu'(X,e_{n+1})=f(X), \ X \in \g
\end{array}
\right.
$$
Then $\g'$ is a solvable Lie algebra, not nilpotent as soon as $f$ is a non-nilpotent derivation.

\begin{center}
\section{The algebraic variety of complex Lie algebras}
\end{center}

\medskip 

A $n-$dimensional complex Lie algebra can be seen as a pair $\frak{g}=(%
\mathbb{C}^{n},\mu )$ where $\mu $ is a Lie algebra law on $\mathbb{C}^{n},$
the underlying vector space to $\frak{g}$ is $\mathbb{C}^{n}$ and $\mu $ the
bracket of $\frak{g}$. We will denote by $L^{n}$ the set of Lie algebra laws
on $\mathbb{C}^{n}$. It is a subset of the vectorial space of alternating
bilinear mappings on $\mathbb{C}^{n}.$
We recall that 
two laws $\mu $ and $\mu ^{\prime }$ $\in L^{n}$ are said isomorphic, if
there is $f\in Gl(n,\mathbb{C)}$ such that
$$
\mu ^{\prime }(X,Y)=f*\mu (X,Y)=f^{-1}(\mu (f(X),f(Y)))
$$
for all $X,Y\in \mathbb{C}^{n}.$ In this case, the Lie algebras $\frak{g}=(\mathbb{C}^{n},\mu )$ and $\frak{g}%
^{\prime }=(\mathbb{C}^{n},\mu ^{\prime })$ are isomorphic.

\bigskip

We have a natural action of the linear group $Gl(n,\C)$ on $L^n$ given by
$$
\begin{array}{cc}
Gl(n,\C)\times  L^n & \longrightarrow L^n \\
( f \ , \mu) &\longrightarrow f*\mu 
\end{array}
$$
We denote by $\mathcal{O(\mu )}$ the orbit of  $\mu$ respect to this action. This orbit is the set of the 
laws isomorphic to $\mu .$

\medskip

\noindent{\bf Notation}. If $\g$ is a $n$-dimensional complex Lie algebra whose bracket (or law) is $\mu $, we 
will denote by $\g=(\mu,\C^n)$ this Lie algebra and often we do not distinguish the Lie algebra and its law. 

\medskip

\subsection{The algebraic variety $L^{n}.$}

Let $\g=(\mu,\C^n)$  a $n$-dimensional complex Lie algebra.  We fix a basis $\{e_{1},e_{2},\cdots ,e_{n}\}$ of $\mathbb{C}^{n}.$ The
structural constants of $\mu \in L^{n}$ are the complex numbers $C_{ij}^{k}$
given by 
$$
\mu (e_{i},e_{j})=\sum_{k=1}^{n}C_{ij}^{k}e_{k}.
$$
As the basis is fixed, we can identify the law $\mu $ with its structural
constants. These constants satisfy :
$$
(1)\left\{
\begin{array}{ll}
C_{ij}^{k}=-C_{ji}^{k}\ ,\ \ \ 1\leq i<j\leq n\ ,\ \ \ 1\leq k\leq n &  \\ 
\\
\sum_{l=1}^{n}C_{ij}^{l}C_{lk}^{s}+C_{jk}^{l}C_{li}^{s}+C_{ki}^{l}C_{jl}^{s}=0\ ,\ \ \ 1\leq i<j<k\leq n\ ,\ \ \ 1\leq s\leq n.
&
\end{array}
\right.
$$
Then $L^{n}$ appears as an algebraic variety embedded in the linear space of
alternating bilinear mapping on $\mathbb{C}^{n}$, isomorphic to $\mathbb{C}^{%
\frac{n^{3}-n^{\acute{2}}}{2}}.$

\medskip

Let be $\mu \in L^{n}$ and consider the Lie subgroup $G_{\mu }$ of $Gl(n,%
\mathbb{C)}$ defined by
$$
G_{\mu }=\{f\in Gl(n,\mathbb{C)}\mid \text{\quad }f*\mu =\mu \}
$$
Its Lie algebra is the Lie algebra of derivations $Der(\frak{g})$ of the Lie algebra $\frak{g}$. We can denote also this algebra by 
$Der(\mu )$ . The orbit $%
\mathcal{O}(\mu )$ is isomorphic to the homogeneous space $Gl(n,\mathbb{C)}$/$G_{\mu
}.\, $Then it is a $\mathcal{C}^{\infty }$ differential manifold of
dimension
$$
\dim \mathcal{O}(\mu )=n^{2}-\dim Der(\mu ).
$$

\medskip

\noindent
\framebox{\parbox{6.25in}{
\begin{proposition}
The orbit $\mathcal{O}(\mu )$ of the law $\mu$ is an homogeneous differential manifold of dimension
$$
\dim \mathcal{O}(\mu )=n^{2}-\dim Der(\mu ).
$$
\end{proposition}
}}

\medskip
\subsection{The tangent space to $\mathcal{O}(\protect\mu )$ at $\protect\mu 
$}

We have seen that the orbit $\mathcal{O}$($\mu $) of $\mu $ is a
differentiable manifold embedded in $L^{n}$ defined by
$$
\mathcal{O}(\mu )=\frac{Gl(n,\mathbb{C)}}{G_{\mu }}
$$
We consider a point $\mu ^{\prime }$ close to $\mu $ in $\mathcal{O}$($\mu $%
). There is $f\in Gl(n,\mathbb{C)}$ such that $\mu ^{\prime }=f*\mu .$
Suppose that $f$ is close to the identity : $f=Id+\varepsilon g,$ with $g\in
gl(n)$. Then
\begin{eqnarray*}
\mu ^{\prime }(X,Y) &=&\mu (X,Y)+\varepsilon [-g(\mu (X,Y))+\mu (g(X),Y)+\mu
(X,g(Y))] \\
&&+\varepsilon ^{2}[\mu (g(X),g(Y))-g(\mu (g(X),Y)+\mu (X,g(Y))-g\mu (X,Y)]
\end{eqnarray*}
and
$$
lim_{\varepsilon \rightarrow 0}\frac{\mu ^{\prime }(X,Y)-\mu (X,Y)}{\varepsilon }=\delta _{\mu }g(X,Y)
$$
where $\delta _{\mu }g$ defined by 
$$\delta _{\mu }g(X,Y)=-g(\mu (X,Y))+\mu (g(X),Y)+\mu(X,g(Y))$$ 
is the coboundary of the cochain $g$ for the Chevalley cohomology of the Lie algebra $\frak{g}$. For simplify the writting, we will denote $B^{2}(\mu ,\mu )$ and $Z^{2}(\mu ,\mu )$ as well as 
$B^{2}(\frak{g} ,\frak{g} )$ and $Z^{2}(\frak{g} ,\frak{g} )$

\medskip

\noindent
\framebox{\parbox{6.25in}{
\begin{proposition}
The tangent space to the orbit $\mathcal{O}$($\mu )$ at the point $\mu $ is
the space $B^{2}(\mu ,\mu )$ of the 2-cocycles of the Chevalley cohomology
of $\mu .$
\end{proposition}
}}

\subsection{The tangent cone to $L^{n}$ at the point $\protect\mu $}

Let $\mu $ be in $L^{n}$ and consider the bilinear alternating mappings $\mu _t=\mu +t\varphi $ 
where $t$ is a complex parameter. Then $\mu _t \in L^{n}$ for all $t$ if and only if we have : 
$$
\left\{ 
\begin{array}{c}
\delta _{\mu}\varphi =0, \\ 
\varphi \in L^{n}.
\end{array}
\right.
$$
So we have the following characterisation of the tangent line of the variety $L^n$:

\medskip 

\noindent
\framebox{\parbox{6.25in}{
\begin{proposition}
A straight line $\Delta $ passing throught $\mu $ is a tangent line in $\mu $
to $L^{n}$ if its direction is given by a vector of $Z^{2}(\mu ,\mu ).$
\end{proposition}
}}

\medskip

Suppose that $H^{2}(\mu ,\mu )=0.$ Then the tangent space to $\mathcal{O}$($%
\mu )$ at the point $\mu $ is the set of the tangent lines to $L^{n}$ at the
point $\mu .$ Thus the tangent space to $L^{n}$ exists in this point and it
is equal to $B^{2}(\mu ,\mu ).$ The point $\mu $ is a nonsingular point. We
deduce of this that the inclusion $\mathcal{O}$($\mu )\hookrightarrow L^{n}$
is a local homeomorphisme. This property is valid for all points of $%
\mathcal{O}$($\mu )$, then $\mathcal{O}$($\mu )$ is open in $L^{n}$ (for the
induced metric topology).

\medskip 

\noindent
\framebox{\parbox{6.25in}{
\begin{proposition}
Let $\mu \in L^{n}$ such that $H^{2}(\mu ,\mu )=0.$ If the algebraic variety
$L^{n}$ is provided with the metric topology induced by $\mathbb{C}^{\frac{%
n^{3}-n^{2}}{2}}$, then the orbit $\mathcal{O}$($\mu )$ is open in $L^{n}.$
\end{proposition}
}}

\medskip

\noindent This geometrical approach shows the significance of problems undelying to the existence of
singular points in the algebraic variety $L^{n}$. We are naturally conduced
to use the algebraic notion of the scheme associated to $L^{n}$. Recall that
the notions of algebraic variety and the corresponding scheme coincide if
this last is reduced.\ In the case of $L^{n}$, we will see that the
corresponding scheme is not reduced.

\bigskip

\begin{center}
\section{The scheme $\mathcal{L}^{n}.$}
\end{center}

\medskip

\subsection{Definition}

Consider the formal variables $C_{jk}^{i}$ with $1\leq i<j\leq n$ and $1\leq
k\leq n$ and let us note $\mathbb{C}[C_{jk}^{i}]$ the ring of polynomials in
the variables $C_{jk}^{i}.$ Let $I_{J}$ the ideal of $\mathbb{C}[C_{jk}^{i}]$
generated by the polynomials associated to the Jacobi relations :
$$
\sum_{l=1}^{n}C_{ij}^{l}C_{lk}^{s}+C_{jk}^{l}C_{li}^{s}+C_{ki}^{l}C_{jl}^{s}
$$
with $1\leq i<j<k\leq n,$ and $1\leq s\leq n.$ The algebraic variety $L^{n}$
is the algebraic set associated to the ideal $I_{J}:$%
$$
L^{n}=V(I_{J})
$$
We will note, for the Jacobi ideal $I_{J}$ of $\mathbb{C}[C_{jk}^{i}]$%
$$
radI_{J}=\{P\in \mathbb{C}[C_{jk}^{i}],\text{ such that }\exists r\in 
\mathbb{N}\text{, }P^{r}\in I_{J}\}
$$
In general, $radI_{J}\neq I_{J}$ (Recall that if $I$ is a maximal ideal,
then $radI=I.$ It is also the case when $I$ is a prime ideal). If $M$ is a
subset of $\mathbb{C}^{\frac{n^{3}-n^{2}}{2}},$ we note by $\frak{i}(M)$ the
ideal of $\mathbb{C}[C_{jk}^{i}]$ defined by
$$
\frak{i}(M)=\{P\in \mathbb{C}[C_{jk}^{i}],P(x)=0\quad \forall x\in M\}.
$$
Then we have
$$
\frak{i}(L^{n})=\frak{i}(V(I_{J}))=radI_{J}
$$
We consider the ring
$$
A(L^{n})=\frac{\mathbb{C}[C_{jk}^{i}]}{I_{J}}
$$
which also is a finite type $\mathbb{C}-$algebra. This algebra corresponds
to the ring of regular functions on $L^{n}.$ Recall that an ideal $I$ of a
ring $A$ is a prime ideal if the quotient ring $A/I$ is an integral ring. In
particular, maximal ideals are prime ideals. The quotient ring is called
reduced if it doesn't contain nonnul nilpotent element ($\forall P\neq
0,\forall n,P^{n}\neq 0)$. As we have generally $radI_{J}\neq I_{J}$, the
algebra $A(L^{n})$ is not reduced.

The affine algebra $\Gamma (L^{n})$ of the algebraic variety $L^{n}$ is the
quotient ring
$$
\Gamma (L^{n})=\frac{\mathbb{C}[C_{jk}^{i}]}{\frak{i}(L^{n})}.
$$
As we have the following inclusion
$$
\frak{i}(L^{n})=rad\frak{i}(L^{n}),
$$
we deduce that the ring $\Gamma (L^{n})$ always is reduced.

Let us note $Spm(A(L^{n}))$ the set of maximal ideals of the algebra $%
A(L^{n}).$ We have a natural bijection between $Spm(A(L^{n}))$ and $L^{n}:$%
$$
Spm(A(L^{n}))\sim L^{n}.
$$
We provide $Spm(A(L^{n}))$ with the Zarisky topology : Let $\frak{a}$ be an
ideal of $A(L^{n})$ and we consider the set $V(\frak{a})$ of maximal ideals
of $A(L^{n})$ containing $\frak{a}$ (in fact we can suppose that $V(\frak{a}%
) $ is the set of radicial ideals containing $\frak{a}$). These sets are the
closed sets of the topology of $Spm(A(L^{n})).$ We  define a basis of
open sets considering, for $f\in A(L^{n})$ :
$$
D(f)=\{x\in Spm(A(L^{n})),\quad f(x)\neq 0\}.
$$
There exists a sheaf of functions $\mathcal{O}_{Spm(A(L^{n}))}$ on $%
Spm(A(L^{n}))$ with values on $\mathbb{C}$ such that, for all $f\in A(L^{n})$,
we have
$$
\Gamma (D(f),\mathcal{O}_{Spm(A(L^{n}))})=A(L^{n})_{f}
$$
where $A(L^{n})_{f}$ is the ring of functions
$$
x\rightarrow \frac{g(x)}{f(x)^{n}}
$$
for $x\in D(f)$ and $g$ $\in A(L^{n}).$ In particular we have
$$
\Gamma (Spm(A(L^{n})),\mathcal{O}_{Spm(A(L^{n}))})=A(L^{n})
$$
The affine scheme of ring $A(L^{n})$ is the space $(Spm(A(L^{n})),\mathcal{O}%
_{Spm(A(L^{n}))})$ noted also $Spm(A(L^{n}))$ or $\mathcal{L}^{n}.$

\subsection{The tangent space of the scheme $\mathcal{L}^{n}$}

The tangent space to the scheme $Spm(A(L^{n}))$ can be calculated
classically. We consider the infinitesimal deformations of the algebra $%
\Gamma (Spm(A(L^{n})),\mathcal{O}_{Spm(A(L^{n}))})$ at a given point. If $%
F_{1},...,F_{N}$ with $N=\frac{1}{6}n(n-1)(n-2)$ are the Jacobi polynomials,
then the tangent space at the point $x$ to the scheme $\mathcal{L}^{n}$ is 
$$
T_{x}(\mathcal{L}^{n})=Ker\,d_{x}(F_{1},..,F_{N})
$$
where $d_{x}$ designates the Jacobian matrix of the $F_{i}$ at the point $x$%
. From the definition of the cohomology of the Lie algebra $\frak{g}$
associated to the point $x,$ we have :

\medskip 

\noindent
\framebox{\parbox{6.25in}{
\begin{theorem}
T$_{x}(Spm(A(L^{n}))=$T$_x(\mathcal{L}^{n})=Z^{2}(\frak{g},\frak{g})$
\end{theorem}
}}

\bigskip 

\begin{center}
\section{Contractions  of Lie algebras}
\end{center}

\medskip
 
Sometimes the word "contraction" is replaced by the word "degeneration". But we prefer here use the term contraction because
its sense is more explicit.

Let us consider the complex algebraic variety $L^{n}$. This variety can be  endowed with the Zariski
topology (the closed set are defined by a finite number of polynomial equations on the parameters $C_{ij}^k$.) It can be also
endowed with the metric topology induced by $\C^N$ where $N=\frac{n^3-n^2}{2}$ the vector space
of structural constants, considering the embeding $L^n \in \C ^N$. Recall that every open set for the Zariski topology is 
an open set for the metric topology. If $A$ is a subset of $L^n$, we note $\bar A$ the closure of $A$ in $L^n$ for the Zariski
topology and $\bar A^d$ its closure for the metric topology.

\medskip

\subsection{Definition of contractions of Lie algebras}
Let $g=(\mu,\C^n)$ be a $n$-dimensional complex Lie algebra.

\medskip 

\noindent
\framebox{\parbox{6.25in}{
\begin{definition}
A Lie algebra $\g _0=(\mu _{0}, \C^n)$, $\mu _0 \in L^{n}$ is called contraction of $\g$ 
if $\mu _{0}\in \overline{\mathcal{O}(\mu )}.$
\end{definition}
}}

\medskip

The historical notion of contraction, given by Segal, was the following. Consider a sequence
$\{f_p\}$ in $Gl(n,\C)$. We deduce a sequence $\{\mu_p\}$ in $L^n$ by putting
$$\mu _p= f_p* \mu.$$
If this sequence admits a limit $\mu _0$ in $\C^N$, then $\mu_0 \in L^n$ and $\mu_0$ was called a contraction of $\mu$. 
The link between these two notions of contraction is given on the following proposition.

\medskip 

\noindent
\framebox{\parbox{6.25in}{
\begin{proposition}
For every $\mu \in L^n$, the Zariski closure $\overline{\mathcal{O}(\mu )}$ of the orbit ${\mathcal{O}(\mu )}$ is equal
to the metric closure $\overline{\mathcal{O}(\mu )}^d$:
$$\overline{\mathcal{O}(\mu )}=\overline{\mathcal{O}(\mu )}^d.$$
\end{proposition}
}}

\medskip

\noindent {\it Proof.} In fact the field of coefficients is $\C$.

\medskip


\noindent
\framebox{\parbox{6.25in}{
\begin{Corollary}
Every contraction of $\mu \in L^n$ is obtained by a Segal contraction.
\end{Corollary}
}}

\medskip

\subsection{ Examples}

\subsubsection{The abelian case} Every Lie algebra can be contracted on the abelian Lie algebra. In fact
if the law $\mu $ is defined on a basis $\{X_{i}\}$ by $\mu
(X_{i},X_{j})=\sum C_{ij}^{k}X_{k},$ we consider the isomorphism $%
f_{\varepsilon }(X_{i})=\varepsilon X_{i},\varepsilon \neq 0.$ Then the law $%
\mu _{\varepsilon }=f_{\varepsilon }\ast \mu $ satisfies $\mu _{\varepsilon
}(X_{i},X_{j})=\sum \varepsilon C_{ij}^{k}X_{k}$ and $lim_{\varepsilon
\rightarrow 0}\mu _{\varepsilon }$ exists and coincides with the law of the
abelian Lie algebra.

\medskip

\subsubsection{Contact Lie algebras} Let us consider the open set $\mathcal{C}_{2p+1}$ of L$^{2p+1\text{ }}$%
constituted of $(2p+1)$-dimensional Lie algebra endowed with a contact form,
that is $\omega \in \frak{g}^{\ast }$ (the dual of $\frak{g}$) satisfying 
$$
\omega \wedge (d\omega )^{p}\neq 0.
$$
There is a basis $(X_{1},X_{2},...,X_{2p+1})$ of $\frak{g}$ such that the
dual basis $(\omega =\omega _{1},\omega _{2}...,\omega _{2p+1})$ satisfies
$$
d\omega _{1}=\omega _{2}\wedge \omega _{3}+\omega _{4}\wedge \omega
_{5}+...+\omega _{2p}\wedge \omega _{2p+1}.
$$
The structural constants respect this basis have the form 
$$
C_{23}^{1}=C_{34}^{1}=...=C_{2p2p+1}^{1}=1.
$$
Consider the isomorphism of $\mathbb{C}^{2p+1}$ given by 
$$
f_{\varepsilon }(X_{1})=\varepsilon ^{2}X_{1},\text{ }f_{\varepsilon
}(X_{i})=\varepsilon X_{i}\quad i=2,...,2p+1.
$$
The structural constant $D_{ij}^{k}$ of $\mu _{\varepsilon }=f_{\varepsilon
}\ast \mu $ respect the basis $\{X_{i}\}$ satisfy 
$$
\left\{ 
\begin{array}{l}
D_{23}^{1}=D_{34}^{1}=...=D_{2p2p+1}^{1}=1 \\ 
D_{ij}^{k}=\varepsilon C_{ij}^{k}\text{ for others indices}
\end{array}
\right.
$$
This implies that $lim_{\varepsilon \rightarrow 0}\mu _{\varepsilon }$
exists and corresponds to the law $\frak{h}_{p}$ of the Heisenberg algebra
of dimension $2p+1.$ Then $\frak{h}_{p}\in \overline{\mathcal{C}_{2p+1}}.$

\medskip 

\noindent
\framebox{\parbox{6.25in}{
\begin{proposition}
Every $2p+1$-dimensional Lie algebra provided with a contact form can be contracted on the $2p+1$-dimensional
Heisenberg algebra $\frak{h}_p$. Moreover, every Lie algebra which is contracted on $\frak{h}_p$ admits a contact form.
\end{proposition}
}}

\medskip

\subsubsection{Frobeniusian Lie algebras} Let $g$ be a $2p$-dimensional Lie algebra. It called frobeniusian if there exists a
non trivial linear form $\omega \in {\frak{g}}^*$ such that $[d\omega]^{p} \neq 0.$ In this case the 2-form $\theta=d\omega$
is an exact symplectic form on $\frak{g}.$ 

\medskip 

\noindent
\framebox{\parbox{6.25in}{
\begin{theorem}
\cite{G1} Let $\{F_{\varphi}\;|\;\varphi  \in\mathbb{C}^{p-1}\}$
be the family on $\left(  p-1\right)  $-parameters of $2p$-dimensional Lie
algebras given by
\[
\left\{
\begin{array}
[c]{l}%
d\omega_{1}=\omega_{1}\wedge\omega_{2}+\sum_{k=1}^{p-1}\omega_{2k+1}%
\wedge\omega_{2k+2}\\
d\omega_{2}=0\\
d\omega_{2k+1}=\varphi_{k}\omega_{2}\wedge\omega_{2k+1},\;1\leq k\leq p-1\\
d\omega_{2k+2}=-\left(  1+\varphi_{k}\right)  \omega_{2}\wedge\omega
_{2k+2},\;1\leq k\leq p-1
\end{array}
\right.
\]
where $\left\{  \omega_{1},..,\omega_{2p}\right\}  $ is a basis of $\left(
\mathbb{C}^{2p}\right)  ^{\ast}$. The family $\{F_{\varphi}\}$ is a complex 
irreducible multiple model for the property ``there exists a linear form whose
differential is symplectic'', that is every $2p$-dimensional frobeniusian complex Lie algebra can be contracted on a Lie algebra belonging
to the family $\{F_{\varphi}\}$.
\end{theorem}
}}

\medskip

It can be easily seen that the algebras $F_{\varphi}$ admit the following graduation: 
if $\left\{X_{1},..,X_{2p}\right\}$ is a dual basis to $\left\{\omega_{1},..,\omega_{2p}\right\}$, 
then $F_{\varphi}= \left(F_{\varphi}\right)_{0}\oplus \left(F_{\varphi}\right)_{1}\oplus \left(F_{\varphi}\right)_{2}$, 
where $\left(F_{\varphi}\right)_{0}=\mathbb{C}{X_{2}}, \left(F_{\varphi}\right)_{1}=\sum_{k=3}^{2p}\mathbb{C}{X_{k}}$ 
and $\left(F_{\varphi}\right)_{2}=\mathbb{C}{X_{1}}$. This decomposition will be of importance for cohomological computations.

\bigskip

\noindent{\it Proof.} Let $\frak g$ be a $2p$-dimensional complex frobeniusian Lie algebra. There exists a basis 
$\{X_1,...,X_{2p}\}$ of $\frak g$ such that the dual basis $\{\omega _1,...,\omega _{2p}\}$ satisfies
$$d\omega _1=\omega _1 \wedge \omega _2 +...+ \omega _{2p-1} \wedge \omega _{2p},$$
i.e. the linear form $\omega _1$ being supposed frobenusian. Let us consider the one parameter change of basis :
$$f_{\epsilon}(X_1)=\epsilon ^2 X_1, \ f_{\epsilon}(X_2)=X_2, \ f_{\epsilon}(X_i)=\epsilon X_i, \ i=3,...,2p.$$
This family defines, when $\epsilon \rightarrow 0$, a contraction ${\frak g}_0$ of $\frak g$ whose Cartan Maurer equations are
$$
\left\{
\begin{array}{l}
d\omega _1 = \omega _1 \wedge \omega _2 +...+ \omega _{2p-1} \wedge \omega _{2p}, \\
d\omega _2 = 0, \\
d\omega _3 = C_{23}^3\omega _2 \wedge \omega _3 +C_{24}^3\omega _2 \wedge \omega _4+...+ 
C_{22p-1}^3\omega _{2} \wedge \omega _{2p-1}+C_{22p}^3\omega _{2} \wedge \omega _{2p}, \\
d\omega _4 = C_{23}^4\omega _2 \wedge \omega _3 +(-1-C_{23}^3)\omega _2 \wedge \omega _4+...+ 
C_{22p-1}^4\omega _{2} \wedge \omega _{2p-1}+C_{22p}^4\omega _{2} \wedge \omega _{2p}, \\
.... \\
d\omega _{2p-1} = C_{22p}^4\omega _2 \wedge \omega _3 +C_{2p}^3\omega _2 \wedge \omega _4+...+ 
C_{22p-1}^{2p-1}\omega _{2} \wedge \omega _{2p-1}+C_{22p}^{2p-1}\omega _{2} \wedge \omega _{2p}, \\
d\omega _{2p} = C_{22p-1}^4\omega _2 \wedge \omega _3 +C_{22p-1}^3\omega _2 \wedge \omega _4+...+ 
C_{22p-1}^{2p}\omega _{2} \wedge \omega _{2p-1}+(-1-C_{22p-1}^{2p-1}\omega _{2} \wedge \omega _{2p}.
\end{array}
\right.
$$
The end of the proof consists to reduce the operator $\psi$ which is defined as the restriction of the adjoint operator
$ad X_2$ to the invariant linear subspace $F$ generated by $\{X_3,...,X_{2p}\}$. We can directely verify the following sentences :

- If $\alpha$ and $\beta$ are eigenvalues of $\psi$ such that $\alpha \neq -1- \beta$, then the eigenspaces $F_{\alpha}$ and
$F_{\beta}$ satisfy $[F_{\alpha},F_{\beta}]=0.$

- If the eigenvalue $\alpha$ of $\psi$ us tot equal to $-\frac{1}{2}$, then, for every $X$ and $Y$ $\in F_{\alpha}$, we have
$[X,Y]=0.$

- If $\alpha$  is an eigenvalue $\alpha$ of $\psi$, then $-1-\alpha$  is also an eigenvalue of $\psi$.

- The multiplicities of the eigenvalues $\alpha$ and $-1-\alpha$ are equal.

- The ordered sequences of Jordan blocks corresponding to the eigenvalues $\alpha$ and $-1-\alpha$ are the same.

\noindent From these remarks, we can find a Jordan basis of $\psi$ such that the matrix of $\psi$ restricted to the invariant
subspace $C_{\alpha} \oplus C_{-1-\alpha}$ where $C_{\lambda}$ designates the characteristic subspaces associated to the eigenvalue
$\lambda$ has the following form :
$$
\left(
\begin{array}{lllllll}
\alpha & 0 & 1 & 0 & 0 & 0 & ...\\
0 & -1-\alpha & 0 & 0 & 0 & 0 & ...\\
0 & 0 & \alpha & 0 & 1 & 0 & ... \\
0 & -1 & 0 & -1-\alpha & 0 & 0 & ... \\
0 & 0 & 0 & 0 & \alpha & 0 & ...\\
0 & 0 & 0 & -1 & 0 & -1 -\alpha & ... \\
\end{array}
\right) .
$$
Then the eigenvalues and their corresponding previous blocks classify the elements of the family of frobeniusian models.

\medskip

\noindent We can find also real models for real frobeniusian algebras.  In this case the previous theorem is written

\medskip 

\noindent
\framebox{\parbox{6.25in}{
\begin{theorem}
\cite{G1} (real case). Let $\{F_{\varphi, \rho}\;|\;\varphi \ , \ \rho \in\mathbb{R}^{p-1}\}$
be the family on $\left(  2p-2\right)  $-parameters of $2p$-dimensional real Lie
algebras given by
$$
\left\{
\begin{array}{l}
\left[ X_1,X_2 \right]=X_1, \\
\left [X_{2k+1},X_{2k+2} \right]= X_1, \ k=0,...,p-1, \\
\left[ X_{2},X_{4k-1} \right]= \varphi _kX_{4k-1} + \rho _kX_{4k+1},\\
\left[ X_{2},X_{4k} \right]= (-1-\varphi _k)X_{4k} - \rho _kX_{4k+2},\\
\left[ X_{2},X_{4k+1} \right]= \rho _kX_{4k-1} + \varphi _kX_{4k+1},\\
\left[ X_{2},X_{4k+2} \right]= \rho _kX_{4k} + (-1-\varphi _k)X_{4k+2},\\
\\
{\mbox{\rm for every}} \ k \leq s , \\
\\
\left[ X_{2},X_{4s+2k-1} \right]= -\frac{1}{2}X_{4s+2k-1} + \rho _{k+s-1}X_{4s+2k},\\
\left[ X_{2},X_{4s+2k} \right]= -\rho _{k+s-1}X_{4s+2k-1} + -\frac{1}{2}X_{4s+2k}, \\
\\
{\mbox{\rm for every}}  \ 2 \leq k \leq p-2  
\end{array}
\right.
$$
where $s$ is a parameter satisfying $ 0 \leq [\frac{p-1}{2}].$
The family $\{F_{\varphi, \rho}\}$ is a real
irreducible multiple model for the property ``there exists a linear form whose
differential is symplectic'', that is every $2p$-dimensional frobeniusian real Lie algebra can be contracted on a Lie algebra belonging
to the family $\{F_{\varphi}\}$.
\end{theorem}
}}

\medskip  

\subsection{In\"{ o}n\"u-Wigner contractions}

The first concept of contractions of Lie algebras has been introduced by Segal, In\"on\"u and Wigner for explain some properties 
related with he classical mechanics, the relativist mechanic and the quantum mechanic. The basic idea is to joint these two last
theories with the classical or the galilean mechanic when the fundamental constants (light velocity, Planck constant) tends to 
infinity or zero.  In this context, a contraction is written $lim_{n \rightarrow \infty}f_n^{-1}[f_n(X),f_n(Y)]$ if this limit 
exists. 
The In\" on\"u-Wigner contractions are a particular case of these limits. We consider a family of isomorphismes 
$\{f_\epsilon\}$ in $Gl(n,\C)$ of the form
$$f_\epsilon = f_1 + \epsilon f_2$$
where $f_1 \in gl(n,\C)$ satisfying is a singular operator $det(f_1)=0$ and $f_2 \in Gl(n,\C)$.  These endomorphims can be easily reduced to the 
following form:
$$
f_1=
\left(
\begin{array}{ll}
Id_r & 0 \\
0 & 0
\end{array}
\right)
\
,
\
f_2=
\left(
\begin{array}{ll}
v & 0 \\
0 & Id_{n-r}
\end{array}
\right)
$$
with $rank(f_1)=rank(v)=r$. 

\noindent Such contractions permits to contracts a given Lie algebra $\frak{g}$ on a Lie algebra $\frak{g} _0$ by 
staying invariant a 
subalgebra $\frak{h}$ of $\frak{g}$ that is $\frak{h}$ is always a subalgebra of $\frak{g} _0$. For example, the homogeneous Lorentz
algebra can be contracted, via an Inonu-Wigner contraction, on the homogeneous Galilean algebra. In the same way, 
the De Sitter
algebras can be contracted on the non-homogeneous Lorentz algebra. Now we give a short description of the Inonu-Wigner contractions.
Let $\g=(\mu ,\C^n)$ be a Lie algebra and let $\frak{h}$ a Lie subalgebra of $\g$. Suppose that $\{e_1,...,e_n\}$
is the fixed basis of $\C^n$ and $\{e_1,...,e_p\}$ is a basis of $\frak{h}$.  Thus 
$$\mu(e_i,e_j)=\sum_{k=1}^p C_{ij}^k e_k, \ \ i,j=1,...,p.$$
Let us consider the In\"on\"u-Wigner isomorphisms given by 
$$
\left\{
\begin{array}{l}
f_\epsilon (e_i)(1+\epsilon )e_i, \ \ i=1,...,p\\
\\
f_\epsilon (e_l)=\epsilon e_l, \ \ l=p+1,...,n.
\end{array}
\right.
$$
Here we have $f_\epsilon =f_1+\epsilon f_2$ with
$$
f_1=
\left(
\begin{array}{ll}
Id_p & 0 \\
0 & 0
\end{array}
\right)
\
,
\
f_2=
\left(
\begin{array}{ll}
Id_p & 0 \\
0 & Id_{n-p}
\end{array}
\right)
.$$
The multiplication $\mu _\epsilon =f_\epsilon *\mu $ is written
$$
\left\{
\begin{array}{l}
\mu _\epsilon (e_i,e_j)=(1+\epsilon )^{-1}\mu (e_i,e_j), \ \ i,j=1,...,p\\
\\
\mu _\epsilon (e_i,e_l)=\epsilon (1+\epsilon )^{-1}\sum_{k=1}^pC_{ij}^k e_k+(1+\epsilon )^{-1}\sum_{k=p+1}^nC_{il}^ke_k, \ \ i=1,..,p, \ l=p+1,...,n\\
\\
\mu _\epsilon (e_l,e_m)=\epsilon^2 (1+\epsilon )^{-1}\sum_{k=1}^pC_{lm}^k e_k+\epsilon \sum_{k=p+1}^nC_{lm}^ke_k, \ \ l,m=p+1,...,n
\end{array}
\right.
$$
If $\epsilon \rightarrow 0$, the sequence $\{\mu_\epsilon \}$ has a limit $\mu_0$ given by
$$
\left\{
\begin{array}{l}
\mu _0 (e_i,e_j)=\mu (e_i,e_j), \ \ i,j=1,...,p\\
\\
\mu _0 (e_i,e_l)=\sum_{k=p+1}^nC_{il}^ke_k, \ \ i=1,..,p, \ l=p+1,...,n\\
\\
\mu _0 (e_l,e_m)=0, \ \ l,m=p+1,...,n.
\end{array}
\right.
$$
The Lie algebra $\g_0=(\mu_0,\C^n)$ is an Inonu-Wigner contraction of $\g$ and $\frak{h}$ is  a subalgebra of $\g_0$. We note
that  the subspace $\C\{e_{p+1},...,e_n\}$ is an abelian subalgebra of $\g_0$. 

\medskip 

\noindent
\framebox{\parbox{6.25in}{
\begin{proposition}
If $\g_0$ is an In\"on\"u-Wigner contraction of $\g$ which led invariant the subalgebra $\frak{h}$ of $\g$ then
$$g_0=\frak{h} \oplus \frak{a}$$
where $\frak{a}$ is an abelian ideal of $\g_0$.
\end{proposition}
}}

\medskip

\noindent{\bf Remarks}

1. If $\frak{h}$ is an ideal of $\g$ then $\g_0=\frak{h}\oplus \frak{a}$ with
$$[\frak{h}, \frak{a}]= 0.$$

\medskip

2. Later Saletan and Levy-Nahas have generalized the notion of In\"on\"u-Wigner contractions conidering 

\noindent a) For Saletan contractions the isomorphismes 
$$f_{\varepsilon}=f_1 + \varepsilon f_2$$
with $det(f_2) \neq 0$. Such an isomorphism can be written
$$f_\varepsilon =\varepsilon Id + (1-\varepsilon)g$$
with $det(g) =0$.  If $q$ is the nilindex  of the nilpotent part of $g$ in its Jordan decomposition, then starting 
with a Lie algebra $\frak{g}$
we contract via $f_\varepsilon $ we obtain a new Lie algebra $\frak{g}_1$, we contract $\frak{g}_1$ via $f_\varepsilon $ and we obtain a Lie algebra 
$\frak{g}_2$ and so on. Thus we construct a sequence of contractions. This sequence is stationnary from the order $q$. 
The In\"on\"u-Wigner case corresponds to $q=1$. 

\noindent b) Levi-Nahas extends the notion of Saletan contractions considering singular contractions. 
In this case the isomorphism $f$ has the following form
$$f=\varepsilon f_1 + (\varepsilon)^2 f_2,$$
$f_1$ and $f_2$ satisfying the Saletan hypothesis.

\medskip

3. In his book, R. Hermann introduce also a notion of contraction but this notion does not correspond to our definition.
In fact he considers some singular contractions where the dimension is not invariant. In our presentation the dimension 
of $\frak{g}$ and its contracted are the same. In the Hermann definition, this is not the case. 

\medskip

4. There exists contraction of Lie algebras which are not In\"on\"u-Wigner contractions. For example 
consider the $4$-dimensional solvable Lie algebra given by 
$$
\left[ e_1,e_2 \right]=e_2, \ \left[ e_3,e_4 \right]=e_4.
$$
This Lie algebra can be contracted on the nilpotent filiform Lie algebra :
$$
\left[ e_1,e_2 \right]=e_3, \ \ \left[ e_1,e_3 \right]=e_4.
$$
From the previous proposition, this contraction cannot be an Inonu-Wigner contraction.

\medskip  

\subsection{In\"on\"u-Wigner contractions of Lie groups}

 There exists a notion of In\"on\"u-Wigner contraction of Lie groups. It is subordinated to the contraction
of its Lie algebras. Then every Lie group can be contracted in the In\"on\"-Wigner sense of its one parameter subgroup. 
The three dimensional rotation group is contracted to the Euclidean group of two dimension. Contraction of the homogeneous 
Lorentz group with respect to the subgroup which leaves invariant the coordinate temporal yields the homogeneous
Galilei group. Contraction of  the inhomogeneous 
Lorentz group with respect to the subgroup generated by the spatial rotations and the time displacements 
yields the full Galilei group. Contraction of the De Sitter group yields the inhomogeneous
Lorentz group. All these exemples are described in the historical paper of In\"on\"u-Wigner \cite{I.N}.

\subsection{Weimar-Woods contractions}

These contractions are given by diagonal isomrphisms
$$f(e_i)=\epsilon ^{n_i}e_i$$
where $n_i \in \Z$ and the contraction is given when $\epsilon \rightarrow 0.$ This  contraction ca be viewed
as generalized In\"on\"u-Wigner contraction with integer exponents and the aim is to construct all the possible contractions.
In the following section we will present the notion of deformation. We will see that any contraction of a Lie algebra $\g$
on $\g_1$ determinates a deformation of $\g_1$ isomorphic to $\g_0$. The notion of Weimer-Woods contractions permits to
solve the reciprocity.

\medskip  

\subsection{The diagrams of contractions}
In the following the symbol $\g_1 \rightarrow \g_2$ means that $\g_1$ can be contracted on $\g_2$ and there is not Lie algebra
$\g$ such that $\g_1$ can be contracted on $\g$ and $\g$ on $\g_3$ (there is no intermediaire).

\noindent The two-dimensional case :

\noindent The variety $L^2$ is the union of $2$ orbits :
$$L^2={\mathcal{O}(\mu _0^2 )}\cup {\mathcal{O}(\mu _1^2 )}$$
where $\mu_0^2$ is the law of the $2$-dimensional abelian Lie algebra and $\mu _1 ^2$
given by 
$$\mu_1^2(e_1,e_2)=e_2.$$
We have 
$$L^2={\overline {\mathcal{O}(\mu _1^2 )}}.$$
Let $\frak{a}_2$ the abalian Lie albegra and $\w_2$ the solvable Lie algebra of law $\mu_1^2$.
The diagramm of contraction is
$$\w_2\longrightarrow \frak{a}_2.$$

\bigskip
\begin{center}
\section{Formal deformations}
\end{center}

\medskip

\subsection{The Gerstenhaber products}

Let $V$ be the $n$-dimensional $\Bbb{C}$-vector space $\Bbb{C}^n$. We denote by $V^{\otimes ^p}$ the tensor product of $V$
with itself $p$ times and $V_p = Hom_{\Bbb{C}}(V^{\otimes ^{p+1}},V)=Bil(V\times V \times...\times V,V)$ the $\Bbb{C}$-vector space
of $(p+1)$-linear forms on $V$ with values on $V$. We have $V_0=End(V)$ and we put $V_{-1}=V$ and $V_{p}=0$ for $p<-1$. We obtain
a sequence $(V_p)_{p \in \Bbb{Z}}$ of  $\Bbb{C}$-vector spaces and 
we note ${\cal{V}}=\oplus _{p \in \Bbb{Z}} V_p$. We define  
products noted $\circ _i$ on $\cal{V}$ putting :

\noindent if $\phi \in V_p, \ \psi \in V_q$ then $\phi \circ _i \psi \in V_{p+q}$ by setting

\medskip 

\noindent
\framebox{\parbox{6.25in}{
$$
\begin{array}{l}
\phi \circ _i \psi (v_0 \otimes ... \otimes v_{i-1} \otimes w_0 \otimes ... \otimes w_q \otimes v_{i+1} \otimes...\otimes v_p) \\
= \phi (v_0 \otimes ... \otimes v_{i-1} \otimes \psi (w_0 \otimes ... \otimes w_q) \otimes v_{i+1} \otimes...\otimes v_p).
\end{array}
$$
}}

\medskip

The system $\{V_m, \circ _i\}$ is a right pre-Lie system, that is we have the following properties where $\phi =\phi ^p$ to
indicate its degree :

$$
(\phi ^p \circ _i \psi ^q) \circ _j \rho ^r =
\left\{
\begin{array}{l}
(\phi ^p \circ _j \rho ^r) \circ _{i+r} \psi ^q \ \ {\mbox {\rm if}} \  \ 0 \leq j \leq i-1 \\
\phi ^p \circ _i (\psi ^q) \circ _{j-i} \rho ^r) \ \ {\mbox {\rm if}}  \ \ i \leq j \leq q+1 ,
\end{array}
\right.
$$
where $\phi , \psi , \rho  \in V_p,V_q,V_r$ respectively. We now define for every $p$ and $q$ a new homomorphism $\circ$ of
$V_p \otimes V_q$ into $V_{p+q}$ by setting for $\phi \in V_p , \psi   \in V_q$

\medskip 

\noindent
\framebox{\parbox{6.25in}{
$$
\phi \circ  \psi =
\left\{
\begin{array}{l}
\phi \circ _0  \psi + \phi \circ _1  \psi +...+ \phi \circ _p  \psi \ \ {\mbox {\rm if}} \ q \ {\mbox {\rm is even}} \\
\phi \circ _0  \psi - \phi \circ _1  \psi +...+(-1)^p \phi \circ _p  \psi \ \ {\mbox {\rm if}} \ q \ {\mbox {\rm is odd}} 
\end{array}
\right.
$$
}}

\medskip

\noindent The vector space ${\cal{V}}=\oplus V_p$ is , endowed with the product $\circ$, a graded pre-Lie algebra. That is, the
graded associator satisfies:

\medskip 

\noindent
\framebox{\parbox{6.25in}{
$$
(\phi ^p \circ \psi ^q) \circ  \rho ^r - (-1)^{pq} (\phi ^p \circ \rho ^r) \circ  \psi ^q =
\phi ^p \circ (\psi ^q \circ  \rho ^r) - (-1)^{pq} \phi ^p \circ (\rho ^r) \circ  \psi ^q).
$$
}}

\medskip

\noindent{ \bf Example} Let be $\mu$ a bilinear mapping on $\Bbb{C}^n$ with values into $\Bbb{C}$. Then 
$$ \mu \circ \mu (X,Y,Z)=\mu (\mu (X,Y),Z) - \mu(X,\mu(Y,Z))$$
and we obtain the associator of the multiplication $\mu.$ This example shows that the class of associative algebras is related by the equation
$$\mu \circ \mu = 0.$$

\medskip

\noindent  {\bf Remark. } 

Starting from the Gerstenhaber products, E.Remm \cite{G.R} intoduces some new products directely related with some non associative 
algebras. Let $\Sigma _n$ be the symmetric group corresponding to the permutations of n elements and $\Sigma _{p+1,q+1}$ with
$p+1+q+1=n$ the subgroup of the $(p+1,q+1)$-shuffles. Let us consider a subgroup $G$ of $\Sigma _{p+1,q+1}$. The Remm product related with $G$
is given by
$$(\phi ^p \circ _G \psi ^q) = \sum _{\sigma \in G} (\phi ^p \circ \psi ^q).\sigma$$
with $\sigma (v_0,...,v_i,w_0,...,w_q,v_{i+1},..,v_p)=(-1)^{\epsilon (\sigma)}(v_{\sigma _1(0)},...,v_{\sigma _1(i)},w_{\sigma _2(0)},...,
w_{\sigma _2(q)},v_{\sigma _1(i+1)},..,v_{\sigma _1(p)}) $ where $\sigma =(\sigma _1, \sigma _2)$ and $\sigma _1 \in 
\Sigma _{p+1}$ and $\sigma _2 \in \Sigma _{q+1}$. 

\noindent The most interesting application concerning these products is in the definition of some non-associative algebras. For
this consider a multiplication on $\Bbb{C}^2$, that is $\mu \in V_1$.  Let $\Sigma _3$ be the symmetric group corresponding to the permutations of three elements. This group contains $6$ subgroups
which are :

. $G_1=\{Id\}$

. $G_2=\{Id,\tau _{12}\}$

. $G_3=\{Id,\tau _{23}\}$

. $G_4=\{Id,\tau _{13}\}$

. $G_5=A_3$ the alternated group

. $G_6=\Sigma _3$

\noindent where $\tau _{ij}$ designates the transposition between $i$ and $j$.For each one of these subgroups 
we define the following non-associative algebra given by the equation:
$$\mu \circ _{G_{i}} \mu =0.$$
Of course, for $i=1$ we obtain an associative multiplication ($\circ _{G_{1}}=\circ $). For $i=2$ we obtain the class of Vinberg 
algebras, for $i=2$ the class of pre-Lie algebras and for $i=6$ the general class of Lie-admissible algebras. The Jacobi condition
corresponds to $i=5$. A large study of these product is made in \cite{R1}, \cite{R2}.

For end this section, we introduce another notation : let be $\phi$ and $\psi$ in $V_1$, and suppose that these bilinear mapping
are alternated. In this case we writte
$$ \varphi \circ \psi = (1/2)\varphi \circ _{G_{5}} \psi.$$ 
Thus
$$
\varphi \circ \psi (X,Y,Z)=\varphi (\psi (X,Y),Z)+\varphi (\psi
(Y,Z),X)+\varphi (\psi (Z,X),Y)
$$
for all $X,Y,Z\in \mathbb{C}^{n}.$ Using this notation, the Lie bracket is
written $\mu \circ \mu =0.$

\medskip

\noindent{\bf Application. }

\noindent Let be $\mu _{0}\in L^{n}$ a Lie algebra law and $\varphi \in C^{2}(\mathbb{C}^{n},\mathbb{C}%
^{n})$ that is a skew-symmetric mapping belonging to $V_1$. Then $\varphi \in Z^{2}(\mu _{0},\mu _{0})$ if and only if
$$
\mu _{0}\circ \varphi +\varphi \circ \mu _{0}=\delta _{\mu _{0}}\varphi =0.
$$

\subsection{Formal deformations}

\medskip

\noindent
\framebox{\parbox{6.25in}{
\begin{definition}
A (formal) deformation of a law $\mu _{0}\in L^{n}$ is a formal sequence
with parameter $t$%
$$
\mu _{t}=\mu _{0}+\sum_{t=1}^{\infty }t^{i}\varphi _{i}
$$
where the $\varphi _{i}$ are skew-symmetric bilinear maps $\mathbb{C}%
^{n}\times \mathbb{C}^{n}\rightarrow \mathbb{C}^{n}$ such that $\mu _{t}$
satisfies the formal Jacobi identity $\mu _{t}\circ \mu _{t}=0.$
\end{definition}
}}

\medskip

Let us develop this last equation.
$$
\mu _{t}\circ \mu _{t}=\mu _{0}\circ \mu _{0}+t\delta _{\mu _{0}}\varphi
_{1}+t^{2}(\varphi _{1}\circ \varphi _{1}+\delta _{\mu _{0}}\varphi
_{2})+t^{3}(\varphi _{1}\circ \varphi _{2}+\varphi _{2}\circ \varphi
_{1}+\delta _{\mu _{0}}\varphi _{3})+...
$$
and the formal equation $\mu _{t}\circ \mu _{t}=0$ is equivalent to the
infinite system 
$$
(I)\text{ }\left\{ 
\begin{array}{l}
\mu _{0}\circ \mu _{0}=0 \\ 
\delta _{\mu _{0}}\varphi _{1}=0 \\
\varphi _{1}\circ \varphi _{1}=-\delta _{\mu _{0}}\varphi _{2} \\ 
\varphi _{1}\circ \varphi _{2}+\varphi _{2}\circ \varphi _{1}=-\delta _{\mu
_{0}}\varphi _{3} \\ 
\vdots \\
\varphi _{p}\circ \varphi _{p}+\sum_{ 1\leq i\leq p-1  }%
\varphi _{i}\circ \varphi _{2p-i}+\varphi _{2p-i}\circ \varphi _{i}=-\delta
_{\mu _{0}}\varphi _{2p} \\ 
\sum_{1\leq i\leq p}\varphi _{i}\circ \varphi _{2p+1-i}+\varphi
_{2p+1-i}\circ \varphi _{i}=-\delta _{\mu _{0}}\varphi _{2p+1} \\ 
\vdots
\end{array}
\right. .
$$
Then the first term $\varphi _{1}$ of a deformation $\mu _{t}$ of a Lie
algebra law $\mu _{0}$ belongs to $Z^{2}(\mu _{0},\mu _{0}).$ This term is
called the infinitesimal part of the deformation $\mu _{t}$ of $\mu _{0}.$

\medskip

\noindent
\framebox{\parbox{6.25in}{
\begin{definition}
A formal deformation of $\mu _{0}$ is called linear deformation if it is of
lenght one, that is of the type $\mu _{0}+t\varphi _{1}$ with $\varphi
_{1}\in Z^{2}(\mu _{0},\mu _{0}).$
\end{definition}
}}

\medskip

For a such deformation we have necessarily $\varphi _{1}\circ \varphi _{1}=0$
that is $\varphi _{1}\in L^{n}.$

\bigskip

\noindent{\bf Examples}

1. Let $\mu _{0}$ be the law of the $n$-dimensional nilpotent Lie algebra defined
by 
$$
\lbrack X_{1},X_{i}]=X_{i+1}
$$
for $i=2,..,n-1$, the other brackets being equal to $0$. Consider a filiform
$n$-dimensional Lie algebra law $\mu $, that is a nilpotent Lie algebra of
which nilindex is equal to $n-1$. Then we prove \cite{G.R} that $\mu $ is isomorphic to an
infinitesimal formal deformation of $\mu _{0}.$ 
\medskip

2. Let us consider a $2p$-dimensional frobeniusian complex Lie algebra. In the previous section
we have done the classification of these algebras up a contraction. In \cite{A.C.1} we prove that
every $2p$-dimensional frobeniusian complex Lie algebra law can be written, up an isomorphism, as
$$\mu = \mu _0 + t \varphi _1$$
where $\mu _0$ is one of a model laws. This prove that every $2p$-dimensional frobeniusian complex Lie algebra is
 a linear deformation of a model.

\bigskip

Now consider $\varphi _{1}\in Z^{2}(\mu _{0},\mu _{0})$ for $\mu _{0}\in
L^{n}.$ It is the infinitesimal part of a formal deformation of $\mu _{0}$
if and only if there are $\varphi _{i}\in C^{2}(\mu _{0},\mu _{0}),$ $i\geq
2 $, such that the system $(I)$ is satisfied. This existence problem is called the {\it formal integration problem}
of $\varphi _1$ at the point $\mu _0$. As the system $(I)$ in an infinite system, we try to solve this by induction.
For $p \geq 2$ let $(I_p)$ the subsystem given by 
$$
(I_p)\text{ }\left\{ 
\begin{array}{l}
\varphi _{1}\circ \varphi _{1}=-\delta _{\mu _{0}}\varphi _{2} \\ 
\varphi _{1}\circ \varphi _{2}+\varphi _{2}\circ \varphi _{1}=-\delta _{\mu
_{0}}\varphi _{3} \\ 
\vdots \\
\sum_{ 1\leq i\leq [p/2]  }%
a_{i,p}(\varphi _{i}\circ \varphi _{p-i}+\varphi _{p-i}\circ \varphi _i)=-\delta_{\mu _{0}}\varphi _{p} \\ 
\end{array}
\right. 
$$
where $a_{i,p}=1$ if $i \neq p/2$ and $a_{i,p}=1/2$ if $i=[p/2]$ (this supposes that $p$ is even). 

\medskip

\noindent
\framebox{\parbox{6.25in}{
\begin{definition}
We say that $\varphi _{1}\in Z^{2}(\mu _{0},\mu _{0})$ is integrable up the order $p$ if there exists 
$\varphi _{i}\in C^{2}(\mu _{0},\mu _{0}),$ $i=2,..p$ such that $(I_p)$ is satisfied.
\end{definition}
}}

\medskip

Suppose that $\varphi _1$ is integrable up the order $p$. Then we prove directly that
$$\sum_{ 1\leq i\leq [p+1/2]  }%
a_{i,p+1}(\varphi _{i}\circ \varphi _{p+1-i}+\varphi _{p+1-i}\circ \varphi _i) \in Z^3(\mu_0,\mu_0).$$
But $\varphi _1$ is integrable up the ordre $p+1$ if and only if this $3$-cochain is in $B^3(\mu_0,\mu_0)$. We deduce the following
result:

\bigskip

\noindent
\framebox{\parbox{6.25in}{
\begin{proposition}
If $H^{3}(\mu _{0},\mu _{0})=0$ then every $\varphi _{1}\in Z^{2}(\mu
_{0},\mu _{0})$ is an infinitesimal part of a formal deformation of $\mu
_{0}.$
\end{proposition}
}}

\medskip

The cohomology class $[$$\sum_{ 1\leq i\leq [p+1/2]  }%
a_{i,p+1}(\varphi _{i}\circ \varphi _{p+1-i}+\varphi _{p+1-i}\circ \varphi _i)]$ is called
the obstruction of index $p+1$. The obstruction of index $2$ are given by the cohomology class of $\varphi _1\circ \varphi _1.$
It can be written using the following quadratic map: 

\medskip

\noindent
\framebox{\parbox{6.25in}{
\begin{definition}
The Rim quadratic map $$sq : H^2(\mu_0,\mu_0)\longrightarrow H^3(\mu_0,\mu_0)$$
is defined by
$$sq([\varphi _1])=[\varphi_1 \circ \varphi _1]$$
for every  $\varphi _{1}\in Z^{2}(\mu _{0},\mu _{0})$. 
\end{definition}
}}

\medskip

\noindent Using this map, the second obstruction is written $sq([\varphi _1])=0.$

\bigskip

\noindent{\bf Remark.} Generalising the notion of formal deformation, we will see in the next section that the infinite system $(I)$ is
equivalent to a finite system, that is there exists only a finite number of obstructions. 

\bigskip

\subsection{Formal equivalence of formal deformations}

Let us consider two formal deformations $\mu _{t}^{1}$ and $\mu _{t}^{2}$ of
a law $\mu _{0}.$ They are called equivalent if there exits a formal linear
isomorphism $\Phi _{t}$ of $\mathbb{C}^{n}$ of the following form 
$$
\Phi _{t}=Id+\sum_{i\geq 1}t^{i}g_{i}
$$
with $g_{i}\in gl(n,\mathbb{C})$ such that 
$$
\mu _{t}^{2}(X,Y)=\Phi _{t}^{-1}(\mu _{t}^{1}(\Phi _{t}(X),\Phi _{t}(Y))
$$
for all $X,Y\in \mathbb{C}^{n}.$

\medskip

\noindent
\framebox{\parbox{6.25in}{
\begin{definition}
A deformation $\mu _{t}^{{}}$ of $\mu _{0}$ is called trivial if it is
equivalent to $\mu _{0}.$
\end{definition}
}}

\bigskip

\noindent In this definition the law $\mu_0$ is considered as a trivial deformation that is when all the bilinear maps
$\varphi _i$ are null.

Let $\mu _{t}^{1}=\mu _{0}+\sum_{t=1}^{\infty }t^{i}\varphi _{i}$ and $\mu
_{t}^{2}=\mu _{0}+\sum_{t=1}^{\infty }t^{i}\psi _{i}$ be two equivalent
deformation of $\mu _{0}.$ It is easy to see that
$$
\varphi _{1}-\psi _{1}\in B^{2}(\mu _{0},\mu _{0}).
$$
Thus we can consider that the set of infinitesimal parts of deformations is
parametrized by $H^{2}(\mu _{0},\mu _{0}).$

Suppose now that $\mu _{t}$ is a formal deformation of $\mu _0$ for which $\varphi _{t}=0$ for
$t=1,..,p$. Then $\delta _{\mu _0}\varphi _{p+1}=0$. If further $\varphi _{p+1} \in  B^{2}(\mu _{0},\mu _{0})$, there exists
$g \in gl(n,\Bbb{C})$ such that $\delta _{\mu _0}g=\varphi _{p+1}$. Consider the formal isomorphism
$\Phi _{t}=Id+tg$. Then 
$$\Phi _{t}^{-1}\mu _t(\Phi _{t},\Phi _{t})=\mu _0 + t^{p+2}\varphi _{t+2} +...$$
and again $\varphi _{t+2} \in  Z^{2}(\mu _{0},\mu _{0})$. 

\medskip

\noindent
\framebox{\parbox{6.25in}{
\begin{theorem}
If $Z^{2}(\mu _{0},\mu _{0})=0$, then $\mu _{0}$ is formally rigid, that is every formal deformation is formally equivalent to
$\mu _0$.
\end{theorem}
}}

\medskip

\noindent {\bf Remarks}

1. In the next chapter we will introduce a notion of topological rigidity, that is the orbit of $\mu _0$ is Zariski open.
The relation between formal deformations and "topological" deformations is done in by the Nijnhuis Richardson theorem. Before to
present this theorem, we will begin to present another notions of deformations, as the perturbations, which are more close of
our topological considerations. The link between Gerstenhaber deformations and perturbations is made in the large context of 
valued deformations (next section).

2. The problem related with the use of formal deformation is the nature of such
tool. We can consider a deformation $\mu _{t}^{{}}$ of a given law $\mu _{0}$
as a $\mathbb{C[[}t\mathbb{]]}$-algebra on $\mathbb{C[[}t\mathbb{]]\otimes }%
_{\mathbb{C}}\mathbb{C}^{n}$ given by
$$
\mu _{t}:\mathbb{C}^{n}\mathbb{\otimes }_{\mathbb{C}}\mathbb{C}%
^{n}\rightarrow \mathbb{C[[}t\mathbb{]]\otimes }_{\mathbb{C}}\mathbb{C}^{n}
$$
where 
$$
\mu _{t}(X\otimes Y)=\mu _{0}(X,Y)+t\otimes \varphi _{1}(X\otimes
Y)+t^{2}\otimes \varphi _{2}(X\otimes Y)+...
$$
But this presentation does not remove the problem concerning the convergence
of the formal serie representing $\mu _{t}.$ The first problem we come up
against is the resolution of the system $(I)$ defined by an infinity of
equations. We will present the deformation differently. Let 
$$
\mu _{t}=\mu _{0}+\sum_{t=1}^{\infty }t^{i}\varphi _{i}
$$
a deformation of $\mu _{0}.$ Amongst the \{$\varphi _{i}\}_{i\geq 1}$ we
extract a free familily $\{\varphi _{i_{1}},\varphi _{i_{2}},\varphi
_{i_{3}},...,\varphi _{i_{N}}\}$ of $C^{2}(\mathbb{C}^{n},\mathbb{C}^{n}).$
We can write again the deformation
$$
\mu _{t}=\mu _{0}+S_{i_{1}}(t)\varphi _{i_{1}}+S_{i_{2}}(t)\varphi
_{i_{2}}+S_{i_{3}}(t)\varphi _{i_{3}}+...+S_{i_{N}}(t)\varphi _{i_{N}}
$$
where $S_{i_{j}}(t)$ is a formal serie of valuation $i_{j}.$ Now it is
evident that such deformation for which one of the formal series $%
S_{i_{j}}(t)$ has a radius of convergence equal to 0 has a limited interest.

\begin{definition}
A convergent formal deformation $\mu _{t}$ of a Lie algebra law $\mu _{0}$ is law of $L^{n}$
defined by 
$$
\mu _{t}=\mu _{0}+S_{1}(t)\varphi _{1}+S_{2}(t)\varphi _{2}+S_{3}(t)\varphi
_{3}+...+S_{p}(t)\varphi _{p}
$$
where

1. $\{\varphi _{1},\varphi _{2},\varphi _{3},...,\varphi _{p}\}$ are
linearly independant in $C^{2}(\mathbb{C}^{n},\mathbb{C}^{n})$

2. The $S_{i}(t)$ are power series with a radius of convergence $r_{i}>0$

3. The valuation $v_{i}$ of $S_{i}(t)$ satisfies $v_{i}<v_{j}$ for $i<j.$
\end{definition}

Let us put $r=\min \{r_{i}\}.$ For all $t$ in the disc $D_{R}=\{t\in \mathbb{%
C},  \mid t\mid \leq R<r\}$ the serie $\mu _{t}$ is uniformly
convergent. With this viewpoint the map
$$
t\in D_{R}\rightarrow \mu _{t}\in L^{n}
$$
appears as an analytic curve in $L^{n}$ passing through $\mu _{0}.$

Let us write the Jacobi conditions concerning $\mu _{t}.$ We obtain
$$
\mu _{t}\circ \mu _{t}=0=S_{1}(t)\delta _{\mu _{0}}\varphi
_{1}+[S_{1}(t)^{2}]\varphi _{1}\circ \varphi _{1}+S_{2}(t)\delta _{\mu
_{0}}\varphi _{2}+...
$$
This identity is true for all $t\in D_{R}.$ By hypothesis, the function
$$
f_{ij}(t)=\frac{S_{j}(t)}{S_{i}(t)t^{j-i}}
$$
is defined by a power serie as soon as $j>i$.\ This function satisfies $%
f_{ij}(0)=0.$ Then

\begin{proposition}
The first term $\varphi _{1}$ of a convergent formal deformation of $\mu _{0}$ satisfies $%
\delta _{\mu _{0}}\varphi _{1}=0.$
\end{proposition}

\begin{definition}
The length of the convergent formal deformation
$$
\mu _{t}=\mu _{0}+S_{1}(t)\varphi _{1}+S_{2}(t)\varphi _{2}+S_{3}(t)\varphi
_{3}+...+S_{p}(t)\varphi _{p}
$$
of $\mu _{0}$ is the integer $p.$
\end{definition}

For example a convergent formal deformation of length $1$ is gievn by
$$
\mu _{t}=\mu _{0}+S_{1}(t)\varphi _{1}.
$$
It is more general that a linear deformation. In fact the deformation
$$
\widetilde{\mu _{t}}=\mu _{0}+\sum_{t=1}^{\infty }t^{i}\varphi _{1}
$$
where all the cochains $\varphi _{i}$ are equal to $\varphi _{1}$
corresponds to the convergent formal deformation of length $1$%
$$
\mu _{t}=\mu _{0}+t\frac{1-t^{n}}{1-t}\varphi _{1}
$$
with $\mid t\mid <1.$

Let us study now the convergent formal deformation of length $2.$ Such deformation can be
writen
$$
\mu _{t}=\mu _{0}+S_{1}(t)\varphi _{1}+S_{2}(t)\varphi _{2}.
$$
The equation of convergent formal deformation
\begin{eqnarray*}
0 &=&S_{1}(t)\delta _{\mu _{0}}\varphi _{1}+S_{1}(t)S_{2}(t)(\varphi
_{1}\circ \varphi _{2}+\varphi _{2}\circ \varphi _{1})+[S_{1}(t)^{2}]\varphi
_{1}\circ \varphi _{1} \\
&&+S_{2}(t)\delta _{\mu _{0}}\varphi _{2}+[S_{2}(t)^{2}]\varphi _{2}\circ
\varphi _{2}
\end{eqnarray*}
gives
$$
\delta _{\mu _{0}}\varphi _{1}=0
$$
and
\begin{eqnarray*}
0 &=&S_{2}(t)(\varphi _{1}\circ \varphi _{2}+\varphi _{2}\circ \varphi
_{1})+[S_{1}(t)]\varphi _{1}\circ \varphi _{1}+\frac{S_{2}(t)}{S_{1}(t)}%
\delta _{\mu _{0}}\varphi _{2} \\
&&+\frac{[S_{2}(t)^{2}]}{S_{1}(t)}\varphi _{2}\circ \varphi _{2}
\end{eqnarray*}
By hypothesis, each fraction is a sum of power series of radius of
convergence equal to $r.$ We can compare the two series $S_{1}(t)$ and $%
\frac{S_{2}(t)}{S_{1}(t)},$ and this depends of the relation between $v_{1}$
and $v_{2}-v_{1}.$ Let us note that $\delta _{\mu _{0}}\varphi _{2}=0$
implies that 
$$
\varphi _{1}\circ \varphi _{2}+\varphi _{2}\circ \varphi _{1}=\varphi
_{1}\circ \varphi _{1}=\varphi _{2}\circ \varphi _{2}=0
$$
and
$$
\mu _{t}=\mu _{0}+S_{1}(t)\varphi _{1}
$$
is a convergent formal deformation of length $1.$ Thus we can suppose that $\delta _{\mu
_{0}}\varphi _{2}\neq 0.$ In this case the 3-cochains $\varphi _{1}\circ
\varphi _{2}+\varphi _{2}\circ \varphi _{1},$ $\varphi _{1}\circ \varphi
_{1},$ $\varphi _{2}\circ \varphi _{2}$ generate a system of rank $1$ in $%
C^{3}(\mathbb{C}^{n},\mathbb{C}^{n}).$ Then the equation of perturbation
implies
$$
\left\{
\begin{array}{l}
\delta _{\mu _{0}}\varphi _{1}=0 \\
\varphi _{1}\circ \varphi _{1}=\delta _{\mu _{0}}\varphi _{2} \\
\varphi _{1}\circ \varphi _{2}+\varphi _{2}\circ \varphi _{1}=\delta _{\mu
_{0}}\varphi _{2} \\ 
\varphi _{2}\circ \varphi _{2}=\delta _{\mu _{0}}\varphi _{2}
\end{array}
\right.
$$

\begin{proposition}
1. Let be $\mu _{t}=\mu _{0}+S_{1}(t)\varphi _{1}$ a convergent formal deformation of length $%
1$.\ Then we have 
$$
\left\{ 
\begin{array}{l}
\delta _{\mu _{0}}\varphi _{1}=0 \\ 
\varphi _{1}\circ \varphi _{1}=0
\end{array}
\right. .
$$
2. If \ $\mu _{t}=\mu _{0}+S_{1}(t)\varphi _{1}+S_{2}(t)\varphi _{2}$ is a
perturbation of length 2 then we have 
$$
\left\{ 
\begin{array}{l}
\delta _{\mu _{0}}\varphi _{1}=0 \\ 
\varphi _{1}\circ \varphi _{1}=\delta _{\mu _{0}}\varphi _{2} \\
\varphi _{1}\circ \varphi _{2}+\varphi _{2}\circ \varphi _{1}=\delta _{\mu
_{0}}\varphi _{2} \\
\varphi _{2}\circ \varphi _{2}=\delta _{\mu _{0}}\varphi _{2}
\end{array}
\right. .
$$
\end{proposition}

\subsection{Perturbations}

The theory of perturbations is based on the infinitesimal framework. The interest of this approach is to have a direct and natural
definition of deformations (which will be called in this context perturbations). We will can consider Lie algebras whose 
structural
constants are infinitelly close to those of a given Lie algebra. This notion can replace the formal notion of deformation and
the convergent formal notion. But it is necessary, for a rigourous description of infinitesimal notions, to work in the context
of Non Standard Analysis. Here we have two ways. The Robinson approach and the Nelson approach. They are different but they have
the same use. In this section we present the notion of perturbation with the Nelson theory because it is more easy. But, in the 
next section, we introducce the notion of valued deformation which generalize the formal deformations and the perturbations. 
But here, the Robinson point of view is better because it is more easy to look the good valuation.  Well we begin to recall what
is the Nelson Non Standard theory.
\subsubsection{Non Standard Theory}
We start with the axiomatic set theory, that is the set theory subordinated with the Zermelo Fraenkel axioms (ZF) more the 
choice axiom (C). Nelson construct a new system of axioms, adding of the previous, three news axioms, noted I,S,T as Idealisation, Sandardisation
and Transfert and a predicat, noted st as standard, in the Zermelo-Fraenkel vocabulary. But the new theory (ZFC + IST) is consistant respect with the ZFC theory. This is the best result of Nelson. It 
permits to do classical mathematic in the IST framework. In the IST theory, the objects construct using only ZFC are called
standard. Then the sets $\Bbb{N}$, $\Bbb{R}$, $\Bbb{C}$ are standard. If we take an element $x \in \Bbb{N}$ which is not standard
then, from the construct of Peano of $\Bbb{N}$, necesary we have $x > n$, for all $n$ standard , $n \in \Bbb{N}$. Such element
will be called infinitely large. We deduce a notion of infinitely large in $\Bbb{R}$ or in $\Bbb{C}$. 

\medskip

\noindent
\framebox{\parbox{6.25in}{
\begin{definition}
Un element $x \in \Bbb{R}$ or in $\Bbb{C}$ is called infinitelly large if it satisfies
$$x > a, \ \forall st(a) \in \Bbb{R} \ \mbox{\rm  or} \ \in \Bbb{C}.$$
It is called infinitesimal if it is zero or if $x^{-1}$ is infinitelly large.

\noindent In other case , it is called appreciable.
\end{definition}
}}

\medskip

These notions  can be extended to $\Bbb{R}^n$ or $\Bbb{C}^n$ for $n$ standard. A vector $v$ of $\Bbb{R}^n$ is infinitesimal
if all its composante are infiniesimal. It is aclled infinitelly large if one of its componants is (this is equivalent to say 
that its euclidean norm is infinitelly large in $\Bbb{R}$).  In the following we denote by IL for an element infinitelly large, 
il for infinitesimal. So we have the following rules

\medskip

\noindent
\framebox{\parbox{6.25in}{
\begin{proposition}

\noindent il + IL = IL.

\noindent IL + IL = IL.

\noindent il $\times$ il = il.

\noindent il / IL = il.
\end{proposition}
}}

\medskip

\noindent On other hand we can in general say nothing about the following cases : il $\times $ IL, IL + IL, IL / IL. 

\medskip

As a consequence of axiom (S) and the completness of $\Bbb{R}$, we have the following important property :

\medskip

\noindent
\framebox{\parbox{6.25in}{
\begin{proposition}
Every appreciable  real element $x$ is infinitely close to a unique standard real element $^ox$, called the shadow of $x$, that is :
$$\forall x \in \Bbb{R}, \ x \ \mbox{\rm appreciable}, \ \exists ! \ ^ox \ \mbox{\rm standard such that}\ x-^ox \ \mbox{\rm inifinitesimal.}$$
\end{proposition}
}}

\medskip

\noindent If $x$ is infinitesimal, its shadow is $0$.

\subsubsection{Perturbations of Lie algebra laws}

Suppose that $\mu _0$ is a standard Lie algebra law on $\Bbb{C}^n$, $n$ being supposed standard.

\medskip

\noindent
\framebox{\parbox{6.25in}{
\begin{definition}
A perturbation $\mu$ of $\mu _0$ is a Lie algebra law on $\Bbb{C}^n$ such that
$$\mu(X,Y) = \mu _0 (X,Y)$$
$\forall X,Y$ standard $\in \Bbb{C}^n$.
\end{definition}
}}

\medskip

Let us consider a standard basis of $\Bbb{C}^n$. Respect this basis, and also respect any standard basis, the structural constants
of $\mu$ and $\mu_0$ are infinitelly close. Then there exists $\epsilon$ infinitesimal such that
$$(\epsilon^{-1}(\mu -\mu _0))(X,Y)$$
are appreciable for all $X,Y$ standard in $\Bbb{C}^n$. Putting 
$$\phi = ^o(\mu -\mu _0),$$this biliniear mapping is skew-symmetric and
we can write
$$\mu = \mu _0 + \epsilon \phi + \epsilon \psi$$
where $\psi$ is a bilinear skew-symmetric mapping satisfying  
$$\psi(X,Y) =0, \ \forall X, Y \ \mbox{ \rm standard} \ \in \Bbb{C}^n. $$
As $\mu$ is a Lie algebra law, it satisfyies the Jacobi condition $\mu \bullet \mu =0.$ Replacing by the expression of $\mu$, 
we obtain
$$(\mu _0 + \epsilon \phi + \epsilon \psi) \bullet (\mu _0 + \epsilon \phi + \epsilon \psi)= 0.$$
This implies, as $\mu _0 \bullet \mu _0 =0$,
$$  \phi \bullet \mu _0 + \mu _0 \bullet \phi + \psi \bullet \mu _0 + \mu _0 \bullet \psi +  \epsilon \phi \bullet \phi 
+ \epsilon \phi \bullet \psi  + \epsilon \psi \bullet \phi +\epsilon \psi \bullet \psi=0.$$
The first part of this equation represents an appreciable trilinear mapping. As it is nul, also its standard part. We deduce that
$$\phi \bullet \mu _0 + \mu _0 \bullet \phi = 2 \delta _{\mu _0} \phi =0.$$

\medskip

\noindent
\framebox{\parbox{6.25in}{
\begin{proposition}
Let $\mu _0$ a standard Lie algebra law on $\Bbb{C}^n$ and let $\mu$ a perturbation of $\mu _0$. Then there exist $\epsilon$
infinitesimal in $\Bbb{C}$ and an infinitesimal bilinear skew-symmetric mapping $\psi$ such that
$$\mu = \mu _0 + \epsilon \phi + \epsilon \psi$$
where $\phi$ is a standardbilinear skew-symmetric mapping satisfying
$$\delta _{\mu _0} \phi =0.$$
\end{proposition}
}}

\medskip

\subsection{Valued deformations}

\subsubsection{Rings of valuation}
We recall briefly the classical notion of ring of valuation. Let $\Bbb{F}$ be a (commutative) field and $A$ a subring of $\Bbb{F}$.
We say that $A$ is a ring of valuation of $\Bbb{F}$ if $A$ is a local integral domain satisfying:
$$ \mbox{ If} \ x \in \Bbb{F} - A, \quad \mbox{then} \quad x^{-1} \in \frak{m}.$$
where $\frak{m}$ is the maximal ideal of $A$.

A ring $A$ is called ring of valuation if  it is a ring of valuation of its field of fractions.

\medskip

\noindent Examples : Let $\Bbb{K}$ be a commutative field of characteristic $0$. The ring of formal series $\Bbb{K}[[t]]$ is a valuation ring. 
On other hand the ring $\Bbb{K}[[t_1,t_2]]$ of two (or more) indeterminates
is not a valuation ring.

\subsubsection{Valued deformations of Lie algebra}
Let $\frak{g}$ be a $\mathbb{K}$-Lie algebra and $A$ a commutative $\mathbb{K}$-algebra of valuation. Then $\frak{g} \otimes A$ is a $\mathbb{K}$-Lie algebra.
We can consider this Lie algebra as an $A$-Lie algebra. We denote this last by $\frak{g}_A$. If $dim_{\mathbb{K}}(\frak{g})$ is finite then
$$dim_{A}(\frak{g}_A)=dim_{\mathbb{K}}(\frak{g}).$$
As the valued ring $A$ is also a  $\mathbb{K}$-algebra  we have a natural embedding of the $\mathbb{K}$-vector space $\frak{g}$ into
the free
$A$-module $\frak{g}_A$. Without loss of generality we can consider this embedding to be the identity map.

\medskip

\noindent
\framebox{\parbox{6.25in}{
\begin{definition}
Let $\frak{g}$ be  a $\mathbb{K}$-Lie algebra and $A$ a commutative $\mathbb{K}$-algebra of valuation such that the residual field
$\frac{A}{\frak{m}}$ is isomorphic to $\mathbb{K}$ (or to a subfield of $\mathbb{K}$). A valued deformation of $\frak{g}$ with base $A$ is a $A$-Lie algebra $\frak{g}'_A$
such that the underlying $A$-module of $\frak{g}'_A$ is $\frak{g}_A$ and that
$$[X,Y]_{\frak{g}'_A} -[X,Y]_{\frak{g}^{ \,} _A} $$ is in the $\frak{m}$-quasi-module $\frak{g} \otimes \frak{m}$ where $\frak{m}$ is the maximal ideal of
$A$.
\end{definition}
}}

\medskip

\noindent{\bf Examples}

1. Formal deformations. The classical notion of deformation studied by Gerstenhaber  is a valued deformation.
In this case $A=\Bbb{K}[[t]]$ and the residual field of $A$ is isomorphic to \K \, . Likewise a versal deformation is a valued deformation.
The algebra $A$ is in this case the finite dimensional \K-vector space $\Bbb{K} \oplus (H^2( \frak{g},\frak{g}))^*$ where $H^2$ denotes the
second Chevalley cohomology
group of $\frak{g}$. The algebra law  is given by
$$(\alpha _1, h_1).(\alpha _2, h_2)=(\alpha _1.\alpha _2, \alpha _1.h_2 + \alpha _2.h_1).$$
It is a local field with maximal ideal $\{0\} \oplus (H^2)^*$. It is also a valuation field because we can endowe this algebra with a field structure,
the inverse of $(\alpha, h)$ being $((\alpha)^{-1}, -(\alpha)^{-2}h)$.

\medskip

2.Versal deformations of Fialowski\cite{F}. Let $\frak{g}$ be a $\mathbb{K}$-Lie algebra and $A$ an unitary commutative local $\mathbb{K}$-algebra. The tensor product 
$ \frak{g} \otimes A $ is naturally
endowed with a Lie algebra structure :
$$[ X \otimes a, Y \otimes b]=[X,Y] \otimes ab.$$
If $ \epsilon : A \longrightarrow \mathbb{K}$, is an unitary augmentation with kernel the maximal ideal $\frak{m}$, a deformation $\lambda $ of $\frak{g}$
 with base $A$ is a Lie algebra structure on $\frak{g} \otimes A$ with bracket $[,]_{\lambda}$ such that
$$ id \otimes \epsilon : \frak{g} \otimes A \longrightarrow
\frak{g} \otimes \Bbb{K}$$
is a Lie algebra homomorphism. In this case the bracket $[,]_{\lambda}$ satisfies
$$[X \otimes 1,Y \otimes 1]_{\lambda}=[X,Y] \otimes 1 + \sum Z_i \otimes a_i$$
where $a_i \in A$ and $X,Y,Z_i \in \frak{g}.$
Such a deformation is called infinitesimal if the maximal ideal $\frak{m}$ satisfies $\frak{m}^2 =0.$ An interesting example is described in [F].
If we consider the commutative algebra $A= \Bbb{K} \oplus
(H^2(\frak{g}, \frak{g}))^*$ (where $^*$ denotes the dual as vector space) such that $dim(H^2) \leq \infty$, the deformation
with base
$A$ is an infinitesimal deformation (which plays the role of an universal deformation).

3. Perturbations. Let $\mathbb{C}^*$ be a non standard extension of $\mathbb{C}$ in the Robinson sense \cite{Ro}. If $\mathbb{C}_l$
is the subring of non-infinitely large elements of  $\mathbb{C}^*$ then the subring $\frak{m}$ of infinitesimals is the 
maximal ideal of
$\mathbb{C}_l$ and $\mathbb{C}_l$ is a valued ring. Let us consider $A=\mathbb{C}_l$. In this case we have a natural embedding of the
variety of $A$-Lie algebras in the variety of $\mathbb{C}$-Lie algebras. Up this embedding (called the transfert principle 
in the Robinson theory), the set of $A$-deformations of $\frak{g}_A$ is an infinitesimal neighbourhood of $\frak{g}$ 
contained in the orbit of $\frak{g}$. Thus any perturbation can be appear as a valued deformation.

\subsection{Decomposition of valued deformations}

In this section we show that every valued deformation can be decomposed in a finite sum (and not as a serie) with pairwise comparable infinitesimal coefficients (that is 
in $\frak{m}$). The interest of this decomposition is to avoid the classical problems of convergence. 
 
\subsubsection{Decomposition in $\frak{m} \times \frak{m}$}
Let $A$ be a valuation ring satisfying the conditions of definition
1. Let us denote by $\cal{F}_A$ the field of fractions of $A$ and
$\frak{m}^2$ the catesian product $\frak{m} \times \frak{m}$ .
Let $(a_1,a_2) \in \frak{m}^2$ with $a_i \neq 0$ for $i=1,2$.

\smallskip

\noindent i)
Suppose that $a_1.a_2^{-1} \in A$. Let be $\alpha = \pi (a_1.a_2^{-1})$ where $\pi$ is the canonical projection on $\frac{A}{\frak{m}}$.
Clearly, there exists a global section
$s: \Bbb{K} \rightarrow A$ which permits to identify $\alpha$ with $s(\alpha)$ in $A$. Then
$$a_1.a_2^{-1}=\alpha + a_3$$
with $a_3 \in \frak{m}$.
Then if $a_3 \neq 0$, 
$$(a_1,a_2)=(a_2(\alpha + a_3),a_2)=a_2(\alpha,1)+a_2a_3(0,1).$$
If $\alpha \neq 0$ we can also write
$$(a_1,a_2)=aV_1+abV_2$$
with $a,b \in \frak{m}$ and $V_1,V_2$ linearly independent in $\Bbb{K}^2$. If $\alpha = 0$ then $a_1.a_2^{-1} \in \frak{m}$ and $a_1=a_2a_3$.
We have
$$(a_1,a_2)=(a_2 a_3,a_2)=ab(1,0)+a(0,1).$$
So in this case, $V_1=(0,1)$ and $V_2=(1,0)$. If $a_3 =0$ then
$$a_1a_2^{-1}=\alpha$$
and 
$$(a_1,a_2)=a_2(\alpha,1)= aV_1.$$

This correspond to the previous decomposition but with $b=0$.

\smallskip

\noindent ii) If
$a_1.a_2^{-1} \in {\cal{F}}_A-A$, then $a_2.a_1^{-1} \in \frak{m}$. We put in this case  $a_2.a_1^{-1}=a_3$ and we have
$$(a_1,a_2)=(a_1,a_1.a_3)=a_1(1,a_3)=a_1(1,0)+a_1a_3(0,1)$$
with $a_3 \in \frak{m}$. Then, in this case the point $(a_1,a_2)$  admits the following decomposition :
$$(a_1,a_2)=aV_1 + ab V_2$$
with $a,b \in \frak{m}$ and $V_1,V_2$ linearly independent in $\Bbb{K}^2$. Note that this case corresponds to the previous but with $\alpha =0$.

\smallskip

Then we have proved

\medskip

\noindent
\framebox{\parbox{6.25in}{
\begin{proposition}
For every point $(a_1,a_2) \in \frak{m}^2$, there exist lineary independent
vectors $V_1$ and $V_2$ in the $\mathbb{K}$-vector space $\mathbb{K} ^2$ such that
$$(a_1,a_2)=aV_1+abV_2$$
for some $a,b \in \frak{m}$. 
\end{proposition}
}}

\medskip

Such decomposition est called of length $2$ if $b \neq 0$. If not it is called of length $1$. 

\subsubsection{Decomposition in $\frak{m}^k$}
Suppose that $A$ is valuation ring satisfying the hypothesis of Definition 1. Arguing as before, we can conclude

\medskip

\noindent
\framebox{\parbox{6.25in}{
\begin{theorem}
For every $(a_1,a_2,...,a_k) \in \frak{m}^k$ there exist $h \ (h \leq k$) independent vectors $V_1,V_2,..,V_h$
whose components are in $\mathbb{K}$ and elements
$b_1,b_2,..,b_h \in \frak{m}$ such that
$$(a_1,a_2,...,a_k)=b_1V_1+b_1b_2V_2+...+b_1b_2...b_hV_h.$$
\end{theorem}
}}

\medskip

The parameter $h$ which appears in this theorem is called the length of the decomposition. This parameter can be different to $k$. It corresponds
to the dimension of the smallest $\mathbb{K}$-vector space $V$ such that $(a_1,a_2,...,a_k) \in V \otimes \frak{m}$. 

\smallskip

\noindent If the coordinates $a_i$ of the vector $(a_1,a_2,...,a_k)$ are in $A$ and not necessarily in its maximal ideal, then writing $a_i=\alpha _i +a'_i$
with $\alpha _i \in \mathbb{K} $ and $a_i' \in \frak{m}$,
we decompose
$$(a_1,a_2,...,a_k)=(\alpha _1,\alpha _2,...,\alpha _k) +(a'_1,a'_2,...,a'_k)$$
and we can apply Theorem 1 to the vector $(a'_1,a'_2,...,a'_k)$.

\subsubsection{Uniqueness of the decomposition}
Let us begin by a technical lemma.

\medskip

\begin{lemma}
Let $V$ and $W$ be two vectors with components in the valuation ring $A$. There exist $V_0$ and $W_0$ with components in $\mathbb{K}$
such that $V=V_0+ V'_0$ and $W=W_0 + W'_0$ and the components of $V'_0$ and $W'_0$ are in the maximal ideal $\frak{m}$. Moreover if
the vectors $V_0$ and $W_0$ are linearly independent then $V$ and $W$ are also independent.
\end{lemma}
{\it Proof.} The decomposition of the two vectors $V$ and $W$ is evident. It remains to prove that the independence of the vectors $V_0$
and $W_0$ implies those of $V$ and $W$. Let $V,W$ be two vectors with components in $A$ such that $\pi (V)=V_0$ and $\pi (W)=W_0$ are independent. Let us suppose that $$xV+yW=0$$ with $x,y \in A$. One of the coefficients $xy^{-1}$ or $yx^{-1}$ is not in $\frak{m}$. Let us suppose that $xy^{-1} \notin \frak{m}$. If $xy^{-1} \notin A$ then $x^{-1}y \in \frak{m}$. Then
$ xV+yW=0$ is equivalent to $V+x^{-1}yW=0$. This implies that $\pi (V)=0$ and this is impossible. Then $xy^{-1} \in A-\frak{m}$. Thus if there exists a linear relation between $V$ and $W$, there exists a linear relation with coefficients in $A-\frak{m}$. We can suppose that $xV+yW=0$ with $x,y \in A-\frak{m}$. As
$V=V_0+V'_0, \, W=W_0+W'_0$ we have
$$\pi (xV+yW)=\pi (x)V_0+\pi (y)W_0=0.$$
Thus $\pi(x)=\pi(y)=0$. This is impossible and the vectors $V$ and $W$ are independent as soon as $V_0$ and $W_0$ are independent vectors. $\quad \Box$

\bigskip

Let $(a_1,a_2,...,a_k)=b_1V_1+b_1b_2V_2+...+b_1b_2...b_hV_h$ and $(a_1,a_2,...,a_k)=c_1W_1+c_1c_2W_2+...+c_1c_2...c_sW_s$ be two decompositions
of the vector $(a_1,a_2,...,a_k)$. Let us compare the coefficients $b_1$ and $c_1$. By hypothesis $b_1c_1^{-1}$ is in $A$ or the inverse is in $\frak{m}$. Then we can suppose that $b_1c_1^{-1} \in A$.
As the residual field is a subfield of \K \, , there exists $\alpha \in \frac{A}{\frak{m}}$ and $c_1 \in \frak{m}$ such that
$$b_1c_1^{-1} = \alpha + b_{11}$$
thus $b_1=\alpha c_1+ b_{11}c_1$. Replacing this term in the decompositions we obtain
$$
\begin{array}{l}
(\alpha c_1+ b_{11}c_1)V_1+(\alpha c_1+ b_{11}c_1)b_2V_2+...+(\alpha c_1+ b_{11}c_1)b_2...b_hV_h \\ =c_1W_1+c_1c_2W_2+...+c_1c_2...c_sW_s.
\end{array}
$$
Simplifying by $c_1$, this expression is written
$$\alpha V_1+ m_1 = W_1 + m_2$$
where $m_1,m_2$ are vectors with coefficients $\in \frak{m}$. From Lemma 1, if $V_1$ and $W_1$ are linearly independent, as its coefficients are in the residual field,
the vectors $\alpha V_1+ m_1$ and $W_1 + m_2$ would be also linearly independent ($\alpha \neq 0$). Thus $W_1=\alpha V_1$. One deduces
$$b_1V_1+b_1b_2V_2+...+b_1b_2...b_hV_h=c_1(\alpha V_1)+c_1b_{11}V_1+c_1b_{12}V_2+...+c_1b_{12}b_3...b_hV_h,$$ with $b_{12}=b_2(\alpha+b_{11}).$
Then
$$b_{11}V_1+b_{12}V_2+...+b_{12}b_3...b_hV_h=c_2W_2+...+c_2...c_sW_s.$$
Continuing this process by induction we deduce the following result
\begin{theorem}
Let be $b_1V_1+b_1b_2V_2+...+b_1b_2...b_hV_h$ and $c_1W_1+c_1c_2W_2+...+c_1c_2...c_sW_s$ two decompositions
of the vector $(a_1,a_2,...,a_k)$. Then

\noindent i. \, $h=s$,

\noindent ii. The flag generated by the ordered free family
$(V_1,V_2,..,V_h)$ is equal to the flag generated by the ordered free
family $(W_1,W_2,...,W_h)$ that is $\forall i \in 1,..,h $
$$ \{
V_1,...,V_i\} =\{ W_1,...,W_i\} $$ where $\{U_i\}$ designates the linear space genrated by the vectors $U_i$.
\end{theorem}

\subsubsection{Geometrical interpretation of this decomposition}
Let $A$ be an $\mathbb{R}$ algebra of valuation. Consider a differential curve $\gamma$ in $\mathbb{R}^3$. We can embed $\gamma$ in a differential curve
$$\Gamma : \mathbb{R} \otimes A \rightarrow \mathbb{R}^3 \otimes A.$$
Let $t=t_0 \otimes 1 + 1 \otimes \epsilon$ an parameter infinitely close to $t_0$, that is $\epsilon \in \frak{m}$. If $M$ corresponds to the point of $\Gamma$ of parameter 
$t$ and $M_0$ those of $t_0$, then the coordinates of the point $M-M_0$ in the affine space $\mathbb{R}^3 \otimes A$ are in $\mathbb{R} \otimes \frak{m}$.
In the flag  associated to the decomposition
of $M-M_0$ we can considere a direct orthonormal frame $(V_1,V_2,V_3)$. It is the Serret-Frenet frame to $\gamma$ at the point $M_0$.

\subsubsection{Decomposition of a valued deformation of a Lie algebra}

Let $\frak{g}'_A$ be a valued deformation with base $A$ of the $\Bbb{K}$-Lie algebra $\frak{g}$. By definition, for every $X$ and $Y$ in $\frak{g}$
we have $[X,Y]_{\frak{g}'_A} -[X,Y]_{\frak{g}^{\,}_A} \in \frak{g} \otimes \frak{m}$. Suppose that $\frak{g}$ is finite dimensional and let $\{X_1,...,X_n \}$
be a basis of $\frak{g}$. In this case
$$[X_i,X_j]_{\frak{g}'_A} -[X_i,X_j]_{\frak{g}^{\,}_A}=\sum_{k} C_{ij}^k X_k$$
with $C_{ij}^k \in \frak{m}$. Using the decomposition of the vector of $\frak{m}^{n^2(n-1)/2}$ with for components $C_{ij}^k$ , we deduce that
$$
\begin{array}{lll}
\lbrack X_i,X_j]_{\frak{g}'_A} -[X_i,X_j]_{\frak{g}^{\,}_A}&=&a\epsilon_1\phi _1(X_i,X_j)+\epsilon_1\epsilon_2\phi _2(X_i,X_j) \\
& &+...+
\epsilon_1\epsilon_2...\epsilon_k\phi _k(X_i,X_j)
\end{array}
$$
where $\epsilon_s \in \frak{m}$ and $\phi _1,...,\phi _l$ are linearly independent. Then we have
$$
\begin{array}{lll}
\lbrack X,Y]_{\frak{g}'_A} -[X,Y]_{\frak{g}^{\,}_A}&=&\epsilon _1\phi _1(X,Y)+\epsilon _1\epsilon _2\phi _2(X,Y) \\
& &+...+
\epsilon _1\epsilon _2...\epsilon _k\phi _k(X,Y)
\end{array}
$$
where the bilinear maps $\epsilon _i$ have values in $\frak{m}$ and linear maps $\phi _i : \frak{g} \otimes \frak{g} \rightarrow \frak{g}$
are linearly independent.

If $\frak{g}$ is infinite dimensional with a countable basis $\{X_n\} _{n \in \mathbb{N}}$ then the $\mathbb{K}$-vector space of
linear map $T_2 ^1= \{ \phi : \frak{g} \otimes \frak{g} \rightarrow \frak{g} \}$ also admits a countable basis. 

\medskip

\noindent
\framebox{\parbox{6.25in}{
\begin{theorem}
If $\mu _{\frak{g}'_A}$ (resp.
$\mu _{\frak{g}^{\,}_A}$) is the  law of the Lie algebra  ${\frak{g}'_A}$ (resp. ${\frak{g}^{\,}_A}$) then
$$\mu _{\frak{g}_A'}-\mu _{\frak{g}^{\,}_A}=\sum _{i \in I} \epsilon _1 \epsilon _2 ...\epsilon _i \phi_i$$
where $I$ is a finite set of indices, $\epsilon _i : \frak{g} \otimes \frak{g} \rightarrow \frak{m}$ are linear maps and $\phi _i$'s are linearly
independent maps in $T_2^1$.
\end{theorem}
}}

\medskip

\subsubsection{Equations of valued deformations}
We will prove that the classical equations of deformation given by Gerstenhaber are still valid in the general frame of valued deformations. Neverless
we can prove that the infinite system described by Gerstenhaber and which gives the conditions to obtain a deformation, can be reduced to a system of finite rank.
Let
$$\mu _{\frak{g}_A'}-\mu _{\frak{g}^{\,}_A}=\sum _{i \in I} \epsilon _1 \epsilon _2 ...\epsilon _i \phi_i$$
be a valued deformation of $\mu$ (the bracket of $\frak{g}$). Then $\mu _{\frak{g}_A'}$ satisfies the Jacobi equations. Following Gerstenhaber we
consider the Chevalley-Eilenberg graded differential complex $\mathcal{C}(\frak{g},\frak{g})$ and the product $\circ$ defined by
$$(g_q \circ f_p)(X_1,...,X_{p+q})=\sum (-1)^{\epsilon (\sigma)} g_q(f_p(X_{\sigma (1)},...,X_{\sigma (p)}),X_{\sigma (p+1)},...,X_{\sigma (q)})$$
where $\sigma$ is a permutation of ${1,...,p+q}$ such that $\sigma (1) < ...<\sigma (p)$ and $\sigma (p+1)<...<\sigma (p+q)$ (it is a 
$(p,q)$-schuffle);  $g_q \in
\mathcal{C}^q(\frak{g},\frak{g})$ and $f_p \in \mathcal{C}^p(\frak{g},\frak{g})$. As $\mu _{\frak{g}_A'}$ satisfies the Jacobi identities,
$\mu _{\frak{g}_A'} \circ \mu _{\frak{g}_A'} = 0$. This gives
$$ (\mu _{\frak{g}_A}+\sum _{i \in I} \epsilon _1 \epsilon _2 ...\epsilon _i \phi_i) \circ ( \mu _{\frak{g}_A}+
\sum _{i \in I} \epsilon _1 \epsilon _2 ...\epsilon _i \phi_i)= 0. \quad \quad (1)$$
As $\mu _{\frak{g}_A} \circ \mu _{\frak{g}_A} = 0$, this equation becomes :
$$\epsilon _1(\mu _{\frak{g}_A} \circ \phi _1 + \phi _1 \circ \mu _{\frak{g}_A})+ \epsilon _1 U=0$$
where $U$ is in $\mathcal{C}^3(\frak{g},\frak{g}) \otimes \frak{m}$. If we symplify by $\epsilon_1$ which is supposed non zero if not the deformation is trivial, we obtain
$$(\mu _{\frak{g}_A} \circ \phi _1 + \phi _1 \circ \mu _{\frak{g}_A})(X,Y,Z)+  U(X,Y,Z)=0 $$
for all $X,Y,Z \in \frak{g}$. As $U(X,Y,Z)$ is in the module $\frak{g} \otimes \frak{m}$ and the first part in $\frak{g} \otimes A$, each one of these
vectors is null. Then $$(\mu _{\frak{g}_A} \circ \phi _1 + \phi _1 \circ \mu _{\frak{g}_A})(X,Y,Z)=0.$$

\medskip

\noindent
\framebox{\parbox{6.25in}{
\begin{proposition}
For every valued deformation with base $A$ of the $\mathbb{K}$-Lie algebra $\frak{g}$, the first term $\phi$ appearing in the associated decomposition is a 2-cochain of the Chevalley-Eilenberg cohomology
of $\frak{g}$ belonging to $Z^2(\frak{g},\frak{g})$.
\end{proposition}
}}

\medskip

We thus rediscover the classical result of Gerstenhaber but in the broader context of valued deformations and not only for the valued deformation
of basis the ring of formal series.

In order to describe the properties of other terms of equations (1) we use the super-bracket of Gerstenhaber which endows the space of
Chevalley-Eilenberg cochains $\mathcal{C}(\frak{g},\frak{g})$ with a Lie superalgebra structure. When $\phi _i \in \mathcal{C}^2(\frak{g},\frak{g}),$
it is defines by
$$[\phi _i,\phi _j]=\phi _i \circ \phi _j + \phi _j \circ \phi _i $$
and $[\phi _i,\phi _j] \in \mathcal{C}^3(\frak{g},\frak{g})$.

\medskip

\noindent
\framebox{\parbox{6.25in}{
\begin{lemma}
Let us suppose that $I=\{1,...,k\}$. If
$$\mu _{\frak{g}_A'}=\mu _{\frak{g}_A}+\sum _{i \in I} \epsilon _1 \epsilon _2 ...\epsilon _i \phi_i$$
is a valued deformation of $\mu$, then the 3-cochains $[\phi _i,\phi _j]$ and $[\mu, \phi_i]$, $1 \leq i,j \leq k-1$, generate 
a linear subspace $V$ of  $\mathcal{C}^3(\frak{g},\frak{g})$
of dimension less or equal to $k(k-1)/2$. Moreover, the 3-cochains $[\phi _i,\phi _j]$, $1 \leq i,j \leq k-1$, form a system of generators of
this space.
\end{lemma}
}}

\medskip
\noindent{\it Proof.}  Let $V$ be the subpace of $\mathcal{C}^3(\frak{g},\frak{g})$ generated by $[\phi _i,\phi _j]$ and $[\mu, \phi_i]$. If $\omega$ is
 a linear form on $V$ of which kernel contains the vectors $[\phi _i,\phi _j]$ for $1 \leq i,j \leq (k-1)$, then the equation (1) gives
$$\epsilon _1 \epsilon _2...\epsilon _k \omega ([\phi _1,\phi _k])+\epsilon _1 \epsilon _2^2...\epsilon _k \omega ([\phi _2,\phi _k])+...+
\epsilon _1 \epsilon _2^2...\epsilon _k^2 \omega ([\phi _k,\phi _k])+\epsilon _2\omega ([\mu,\phi _2]) $$
$$ +\epsilon _2\epsilon _3 \omega ([\mu,\phi _3])...+\epsilon _2\epsilon _3...\epsilon _k \omega ([\mu,\phi _k])=0.$$
As the coefficients which appear in this equation are each one in one $\frak{m}^p$, we have necessarily
$$\omega ([\phi _1,\phi _k])=...=\omega ([\phi _k,\phi _k])=\omega ([\mu,\phi _2])=...=\omega ([\mu,\phi _k])=0$$
and this for every linear form $\omega$ of which kernel contains $V$. This proves the lemma.

From this lemma and using the descending sequence
$$\frak{m} \supset \frak{m}^{(2)} \supset ... \supset \frak{m}^{(p)} ...$$
where $\frak{m}^{(p)}$ is the ideal generated by the products ${ a_1a_2...a_p, \  a_i \in \frak{m} }$ of length $p$,  we obtain :

\medskip

\noindent
\framebox{\parbox{6.25in}{
\begin{proposition}
If
$$\mu _{\frak{g}_A'}=\mu _{\frak{g}_A}+\sum _{i \in I} \epsilon _1 \epsilon _2 ...\epsilon _i \phi_i$$
is a valued deformation of $\mu$, then we have the following linear system :
$$
\left\{
\begin{array}{l}
\delta \phi _1 = 0 \\
\delta \phi _2 = a_{11}^2 [\phi _1,\phi _1] \\
\delta \phi _3 = a_{12}^3 [\phi _1,\phi _2]+a_{22}^3[\phi _1,\phi _1]\\
... \\
\delta \phi _k = \sum _{1 \leq i \leq j \leq k-1} a_{ij}^{k} [\phi _i,\phi _j] \\
\lbrack \phi _1,\phi _k]=\sum _{1 \leq i \leq j \leq k-1} b_{ij}^{1} [\phi _i,\phi _j] \\
.... \\
\lbrack \phi _{k-1},\phi _k]=\sum _{1 \leq i \leq j \leq k-1} b_{ij}^{k-1} [\phi _i,\phi _j]

\end{array}
\right.
$$
where $\delta \phi _i = [\mu,\phi _i]$ is the coboundary operator of the Chevalley cohomology of the Lie algebra $\frak{g}$.
\end{proposition}
}}

\medskip

Let us suppose that the dimension of $V$ is the maximum $k(k-1)/2$. In this case we have no other relations between the generators
of $V$ and the previous linear system is complete, that is the equation of deformations does not give other relations than the relations of this system. 
The following result shows that, in this case, such deformation is
isomorphic, as Lie algebra laws, to a "polynomial"
valued deformation.

\medskip

\noindent
\framebox{\parbox{6.25in}{
\begin{proposition}
Let be $\mu _{\frak{g}_A'}$ a valued deformation of $\mu$ such that
$$\mu _{\frak{g}_A'}=\mu _{\frak{g}_A}+\sum _{i=1,...,k} \epsilon _1 \epsilon _2 ...\epsilon _i \phi_i$$
and dim$V$=k(k-1)/2. Then there exists an automorphism of $\mathbb{K}^n \otimes \frak{m}$ of the form $f=Id \otimes P_k (\epsilon)$ with
$P_k(X) \in \mathbb{K}^{k}[X]$ satisfying $P_k(0)=1$ and $\epsilon \in \frak{m}$ such that the valued deformation $\mu _{\frak{g}''_A}$ defined by
$$\mu _{\frak{g}''_A}(X,Y)=f^{-1}(\mu _{\frak{g}'_A}(f(X),f(Y)))$$
is of the form
$$\mu _{\frak{g}_A"}=\mu _{\frak{g}_A}+\sum _{i=1,...,k} \epsilon ^i \varphi_i$$
where $\varphi _i = \sum _{j \leq i} \phi _j.$
\end{proposition}
}}

\medskip

\noindent {\it Proof.
 } Considering the Jacobi equation
$$ \lbrack\mu _{\frak{g}_A'},\mu _{\frak{g}_A'}]=0$$
and writting that dim$V$=$k(k-1)/2$, we deduce that there exist polynomials $P_i(X) \in \mathbb{K}[X]$ of degree $i$ such that
$$\epsilon _i = a_i \epsilon _k \frac{P_{k-i}(\epsilon _k)}{P_{k-i+1}(\epsilon _k)}$$ with $a_i \in \mathbb{K}$.
Then we have
$$\mu _{\frak{g}_A'}=\mu _{\frak{g}_A}+\sum _{i=1,...,k} a_1a_2...a_i (\epsilon _k)^i \frac{P_{k-i}(\epsilon _k)}{P_{k}(\epsilon _k)} \phi_i.$$
Thus
$$ P_{k}(\epsilon _k)\mu _{\frak{g}_A'}=P_{k}(\epsilon _k)\mu _{\frak{g}_A}+\sum _{i=1,...,k} a_1a_2...a_i (\epsilon _k)^i P_{k-i}(\epsilon _k) \phi_i.$$
If we write this expression according the increasing powers we obtain the announced expression. $\quad \Box$

\noindent Let us note that, for such deformation we have
$$
\left\{
\begin{array}{l}

\delta \varphi _2 + [\varphi _1,\varphi _1] =0\\
\delta \varphi _3 + [\varphi _1,\varphi _2]=0\\
... \\
\delta \varphi _k + \sum _{i+j = k}  [\varphi _i,\varphi _j] =0\\
\sum _{i+j=k+s}  [\varphi _i,\varphi _j] =0.\\

\end{array}
\right.
$$

\subsubsection{Particular case : one-parameter deformations of Lie algebras}
In this section the valuation ring $A$ is $\Bbb{K}[[t]]$. Its maximal ideal is $t\Bbb{K}[[t]]$ and the residual field is $\Bbb{K}$. 
Let $\frak{g}$ be
a $\mathbb{K}$- Lie algebra. Consider $ \frak{g} \otimes A$ as an $A$-algebra and let be $\frak{g}_A'$ a valued deformation of $\frak{g}$.
The bracket $[,]_t$ of this Lie algebra satisfies
$$[X,Y]_t=[X,Y] + \sum t^i \phi _i (X,Y).$$
Considered as a valued deformation with base $\mathbb{K}[[t]]$, this bracket can be written
$$[X,Y]_t = [X, Y] + \sum _{i=1}^{i=k}c_1(t)...c_i(t) \psi _i (X,Y)$$
where $(\psi _1,...,\psi _k)$ are linearly independent and $c_i(t) \in t\mathbb{C}[[t]]$. As $\phi _1 = \psi _1$, this bilinear map belongs
to $Z^2(\frak{g}, \frak{g})$ and we find again the classical result of Gerstenhaber. Let $V$ be the $\mathbb{K}$-vector space generated by
 $[\phi _i,\phi _j]$ and $[\mu, \phi _i], \ i,j=1,...,k-1$, $\mu$ being the law of $\frak{g}$. If dim$V =
k(k-1)/2$ we will say that one-parameter deformation $[,]_t$ is of maximal rank.

\medskip

\noindent
\framebox{\parbox{6.25in}{
\begin{proposition}
Let
$$[X,Y]_t=[X,Y] + \sum t^i \phi _i (X,Y)$$
be a one-parameter deformation of $\frak{g}$. If its rank is maximal then this deformation is equivalent to a polynomial deformation
$$[X,Y]_t'=[X,Y]  + \sum _{i=1,...,k}t^i \varphi _i$$
with $\varphi _i =\sum _{j=1,...,i}a_{ij} \psi _j.$
\end{proposition}
}}

\bigskip

\medskip

\noindent
\framebox{\parbox{6.25in}{
\begin{Corollary}
Every one-parameter deformation of maximal rank is equivalent to a local non valued deformation with base the local algebra $\mathbb{K}[t]$.
\end{Corollary}
}}

\medskip

Recall  that the algebra $\mathbb{K}[t]$ is not an algebra of valuation. But every local ring is dominated by a valuation ring. Then this corollary can be interpreted as saying that every deformation in the local
algebra $\mathbb{C}[t]$ of polynomials with coefficients in $\mathbb{C}$ is equivalent to a "classical"-Gerstenhaber deformation with maximal rank.

\section{The scheme $\mathcal{L}^{n}$ and the cohomology $H(L^{n},L^{n})$}

In this section, we are interested by geometrical properties of $\mathcal{L}%
^{n}$ and $L^{n}$ which can be expressed using the Chevalley cohomology $H(%
\frak{g},\frak{g}).$ Let $x$ be a point of the scheme $\mathcal{L}^{n}$ and $%
\frak{g\in }L^{n}$ the corresponding Lie algebra. The algebraic group $Gl(n,%
\mathbb{C})$ operates on $L^{n}$ and we have noted by $\mathcal{O(}\frak{g}$%
) the orbit of $\frak{g}$. The scheme corresponding to the algebraic variety 
$Gl(n,\mathbb{C})$ is reduced (every group scheme is reduced because the
characteristic of the field is $0$) and is nothing that as the variety
itself. Then, if we denote by $T_{x}^{0}(\mathcal{L}^{n})$ the tangent space
to the orbit scheme of $x$ at the point $x$, we have :

\medskip

\noindent
\framebox{\parbox{6.25in}{
\begin{proposition} \cite{Ga}
There is a canonical isomorphism
$$
\frac{T_{x}(\mathcal{L}^{n})}{T_{x}^{0}(\mathcal{L}^{n})}\simeq H^{2}(\frak{g%
},\frak{g}).
$$
\end{proposition}
}}

\medskip

Let $x$ be in $\mathcal{L}^{n}$ and $\frak{g}$ the corresponding Lie algebra.

\begin{lemma}
If H$^{p}(\frak{g},\frak{g})=0,$ then there exits a neigborhood $U$ of $x$
in $\mathcal{L}^{n}$ such as H$^{p}(\frak{g}_{y},\frak{g}_{y})=0$ for all
Lie algebra $\frak{g}_{y}$ corresponding to the points $y\in U.$
\end{lemma}

Indeed, the function 
$$
\delta ^{i}:x\in \mathcal{L}^{n}\rightarrow \dim H^{i}(\frak{g},\frak{g})
$$
is upper semicontinuous. Then, if H$^{p}(\frak{g},\frak{g})=0,$ this implies
that there is a neigborhood of $x$ suc that $dimH^{p}(\frak{g}_{y},\frak{g}%
_{y})=0$ for every Lie algebra $\frak{g}_{y}$ associated to the points $y$
of $U$.

\medskip

\noindent
\framebox{\parbox{6.25in}{
\begin{theorem}
Let $\frak{g}$ be in $L^{n}$ with $H^{1}(\frak{g},\frak{g})=0.$ Then there
is an open neigborhood $U$ of $x$ in $\mathcal{L}^{n}$ such as the dimension
of the orbits of the points $y\in U$ (and then of the Lie algebras $\frak{g}%
_{y}$) are constant and equal to $n^{2}-n+\dim H^{0}(\frak{g},\frak{g}).$
\end{theorem}
}}

\medskip

Indeed, if $H^{1}(\frak{g},\frak{g})=0,$ from the previous lemma, $H^{1}(%
\frak{g}_{y},\frak{g}_{y})=0$ for every Lie algebra $\frak{g}_{y}$
associated to the points $y$ of $U$. Thus we have 
$$
\dim H^{0}(\frak{g}_{y},\frak{g}_{y})=\dim H^{0}(\frak{g},\frak{g})
$$
As the codimension of the orbit of $\frak{g}_{y}$is $n-\dim $ $H^{0}(\frak{g}%
_{y},\frak{g}_{y})$ , we have $codim\mathcal{O(}\frak{g}_{y})=co\dim
\mathcal{O(}\frak{g})=n-\dim $ $H^{0}(\frak{g},\frak{g}).$

\medskip

\noindent
\framebox{\parbox{6.25in}{
\begin{Corollary}
If $\frak{g}$ satisfies $H^{1}(\frak{g},\frak{g})=0$ and if $\frak{g}$ is a
contraction of $\frak{g}^{\prime }$, then $\frak{g}$ and $\frak{g}^{\prime }$
are isomorphic.
\end{Corollary}
}}

\medskip

This follows from the definition of the contraction : $\frak{g\in }\overline{%
\mathcal{O(}\frak{g}^{\prime })}.$ As $U\cap \mathcal{O(}\frak{g}^{\prime
})\neq \emptyset ,$ $\dim \mathcal{O(}\frak{g}^{\prime })=\dim \mathcal{O(}%
\frak{g}).$ Thus $\mathcal{O(}\frak{g}^{\prime })=\mathcal{O(}\frak{g}).$

\bigskip

\begin{center}
\section{Rigid Lie algebras}
\end{center}

\bigskip

\subsection{Definition}

Let $\g=(\C ^n,\mu) $ be a $n$-dimensional complex Lie algebra. We note always by $\mu $ the corresponding point of $L_n$.

\medskip

\noindent
\framebox{\parbox{6.25in}{
\begin{definition}
The Lie algebra $\g$ is called rigid if its orbit $\mathcal{O}(\mu )$ is open (for the Zariski topology) in $L_n$.
\end{definition}
}}

\medskip

\noindent{\bf Example}. Let us consider the $2$-dimensional Lie algebra given by 
$$\mu(e_1,e_2)=e_2.$$
We have seen in the previous section that its orbit is the set of Lie algebra products given by
$$\mu_{a,b}(e_1,e_2)=ae_1+be_2$$
for all $a,b \in \C.$ This orbit is open then $\mu$ is rigid.

\medskip

\noindent{\bf Remark.} It is no easy to work with the Zariski topology (the closed set are given by polynomial equations).
But, in our case, we can bypass the problem. Indeed, the orbit of a point is an homogeneous space and thus provided with 
a differential structure. In this case the topoly of Zariski and the metric topology induced by the metric topology 
of the vector space $\C^N$ of structure constants coincide. Then $\mu$ is rigid if $\mathcal{O}(\mu )$ is open
in $\C^N$. 

\medskip

\subsection{Rigidity and Deformations}

Intuitively a Lie algebra is rigid if any close Lie algebra is the same or isomorphic. We have presented the
notion of deformations as a notion which try to describe what is a close Lie algebra.  Then we shall have criterium of rigidity
in terms of deformations or perturbations. We recall that a perturbation is a special valued deformation whose
valued ring of coefficients is given by the ring of limited elements in a non standard extension of the complex nuber.

\medskip

\noindent
\framebox{\parbox{6.25in}{
\begin{theorem}
Let $\g=(\C^n,\mu)$ a $n$-dimensional complex Lie algebra. Then $\g$ is rigid if and only if any perturbation
$\g'$ ( where the structure constants are in the ring of limited elements of the Robinson extension $\C ^*$ of $\C$)
is ($\C^*$)-isomorphic to $\g$.
\end{theorem}
}} 

\medskip

\noindent{\bf Example.} Let us consider the Lie algebra $sl(2)$. It is the $3$-dimensional Lie algebra given by
$$
\left\{
\begin{array}{l}
\mu(e_1,e_2)=2e_2\\
\mu(e_1,e_3)=-2e_3\\
\mu(e_2,e_3)=e_1.
\end{array}
\right.
$$
This Lie algebra is simple. Let us prove that it is rigid. Let $\mu'$ be a perturbation of $\mu$. Then the linear 
operator
$ad_{\mu'}e_1$ (of the vector space $(\C^*)^3$ is a perturbation of $ad_{\mu'}e_1$. Then the eigenvalue are close. This implies
that the eigenvalue of $ad_{\mu'}e_1$ are $\epsilon _0,2+\epsilon _1,-2+\epsilon _2$ where $\epsilon _i$ are in the maximal
ideal, that is are infinitesimal. As $ad_{\mu'}e_1(e_1)=0$, we have $\epsilon _0=0$. There exists a basis 
$\{e_1,e'_2,e'_3\}$ of $(\C^*)^3$ whose elements are eigenvectors of $ad_{\mu'}e_1$. We can always choose $e'_2$ such that
$e'_2-e_2$ has coefficients in the maximal ideal. Likewise for $e'_3$. Thus $\mu'(e'_2,e'_3)$ is an eigenvector associated
to $2+\epsilon _1-2+\epsilon _2=\epsilon _1+\epsilon _2$. But this vector is linealrly dependant with $e_1$. Then 
$\epsilon _1+\epsilon _2=0$. We have then
$$
\left\{
\begin{array}{l}
\mu'(e_1,e'_2)=2(1+\epsilon _1)e'_2\\
\mu'(e_1,e'_3)=-2(1+\epsilon _1)e'_3\\
\mu'(e'_2,e'_3)=e_1.
\end{array}
\right.
$$
If we put $e_1=(1+\epsilon _1)e'_1$ we obtain
$$
\left\{
\begin{array}{l}
\mu'(e_1,e'_2)=2e'_2\\
\mu'(e_1,e'_3)=-2(1+e'_3\\
\mu'(e'_2,e'_3)=e_1.
\end{array}
\right.
$$
and $\mu'$ is isomorphic to $\mu$. The Lie algebra $sl(2) $ is rigid.

We similar arguments, we can generalize this result :
\medskip

\noindent
\framebox{\parbox{6.25in}{
\begin{theorem}
Every complex simple or semi-simple Lie algebra is rigid.
\end{theorem}
}}

\medskip

\noindent The proof can be read in \cite{L.G}.

\medskip

\subsection{Structure of rigid Lie algebras}

\medskip

\medskip

\noindent
\framebox{\parbox{6.25in}{
\begin{definition}
A complex (or real) Lie algebra $\g$ is called algebraic if $\g$ is isomorphic
to a Lie algebra of a linear algebraic Lie group.
\end{definition}
}}

\medskip

A linear algebraic Lie group is a Lie subgroup of $GL(p,\C)$ which is defined as the set of zeros of a finite system
of prolynomial equations. Contrary to the variety $L_n$, none of the points defined by this system is singular. This implies that
an algebraic Lie group is a differential manifold. 

\medskip

\noindent{\bf Examples.}

1. Every simple Lie algebra is algebraic.

2. Every nilpotent Lie algebra, or every Lie algebra whose radical is nilpotent is algebraic.

3. Every $n$-dimensional Lie algebra whose the Lie algebra of Derivations is also of dimension $n$ is algebraic.

4. Every Lie algebra satisfying $\mathcal{D}^1(\g)=\g$ is algebraic.

\medskip

We can see that the family of algebraic Lie algebras, for a given dimension, is very big and we have 
not general classification.

\medskip

\noindent
\framebox{\parbox{6.25in}{
\begin{proposition}
The following properties are equivalents :

1. $\g$ is algebraic.

2. $ad(\g)=\{ad_\mu X, X \in \g\}$ is algebraic

3. $\g=\s \oplus \n \oplus \m$ where $\s$ is the Levi semi-simple part, $\n$ the nilradical,that is the maximal nilpotent ideal
and $t$ a Malcev torus such that $ad_\mu \m$ is algebraic.
\end{proposition}
}}

\medskip

In this proposition we have introduced the notion of Malcev torus. It is an abelian subalgebra $\m$ of $ \g$ such that
the endomorphisms  $adX$, $X \in \m$ are semi-simple (simultaneously diagonalizable). All the maximal
torus for the inclusion are conjugated and their commun dimension is called the rank.

\medskip

\noindent
\framebox{\parbox{6.25in}{
\begin{theorem}
Every rigid complex Lie algebra is algebraic.
\end{theorem}
}}

\medskip

This theorem, proposed by R. Carles in \cite{Ca} has many important consequence. For example, we  prove in \cite{A.G} that the
the torus $\m$ is maximal. This permits to construct many rigid Lie algebra and to present some classification. For
example in \cite{G.A} we give the classification of rigid solvable Lie algebras whose nilradical is filiform. We
give also the general classification of $8$-dimensional rigid Lie algebras.

\medskip

\subsection{Rigidity and cohomology}

This approach is more geometrical. If $\g$ is rigid, then its orbit is open. This implies that the tangent spaces to the orbit
at $\mu$ and to ths schema $\mathcal{L}_n$ always at $\mu$ coincide. Before we have determinated these vector spaces. 
We saw that they coincided with spaces $B^2(\g,\g)$ and $Z^2(\g,\g)$. We deduce the following Nijenhuis-Richardson theorem:

\medskip

\noindent
\framebox{\parbox{6.25in}{
\begin{theorem}
If $H^2(\g,\g)=0$ then $\g$ is rigid.
\end{theorem}
}}

\medskip

For example any semi simple Lie algebra has a trivial space of cohomology. Then it is rigid. Every rigid of dimension
less or equal to $8$ satisfies $H^2(\g,\g)=0.$ But there exists many examples of rigid Lie algebras which satisfy
$H^2(\g,\g) \neq 0.$ (see for example \cite{G.A}). A geometric characterization of this fact is given in the following:

\medskip

\noindent
\framebox{\parbox{6.25in}{
\begin{proposition}
Let $\g$ be a rigid complex Lie algebra such that $H^2(\g,\g) \neq \{0\}.$ Then the schema $\mathcal{L}_n$ is not reduced
at the point $\mu$ corresponding to $\g$.
\end{proposition} 
}}

\bigskip

\noindent{\bf Remark.} This last time a study of real rigid Lie algebra has been proposed in \cite{A.C}.

\end{document}